\title{Continuous Hochschild Cohomology and Formality}
\author{
	Patrick Antweiler\thanks{Email: patrick.antweiler"AT"uni-hamburg.de} \\
	Faculty of Mathematics, Universitity of Hamburg
}
\documentclass{article}
\usepackage[top=2.5cm,bottom=2.5cm,left=4cm,right=4cm]{geometry}
\usepackage{amsfonts}

\usepackage{graphicx}
\usepackage{float}
\usepackage{amsmath}
\usepackage{amssymb}
\usepackage{amsthm}
\usepackage{bbm}
\usepackage{hyperref}
\usepackage{enumitem} 
\usepackage[utf8]{inputenc}
\usepackage[T1]{fontenc}
\usepackage{lmodern}
\usepackage{mathrsfs}
\usepackage{changepage}
\usepackage{tikz-cd} 
\usepackage{multirow}
\usepackage[mathscr]{eucal}
\usepackage[toc,page]{appendix}
\usepackage{wasysym}
\usepackage{adjustbox}
\usepackage{blindtext}
\usepackage{makecell}
\usepackage{url}
\usepackage{pgfplots}
\usepackage{tikz}
\usepackage{tikz-3dplot}
\pgfplotsset{compat=1.15}
\usepackage{mathrsfs}
\usetikzlibrary{arrows}

\usepackage{titling}
\usepackage{titling}

\newenvironment{nalign}{
	\begin{equation}
		\begin{aligned}
		}{
		\end{aligned}
	\end{equation}
	\ignorespacesafterend
}

\theoremstyle{definition}

\newtheorem{def1}{Definition}[subsection]

\newtheorem{remark}[def1]{Remark}
\newtheorem{example}[def1]{Example}
\theoremstyle{theorem}
\newtheorem{lemma}[def1]{Lemma}
\newtheorem{corollary}[def1]{Corollary}
\newtheorem{theorem}[def1]{Theorem}
\newtheorem{theorem*}{Theorem}
\newtheorem{proposition*}{Proposition}
\newtheorem*{theorem**}{Theorem}
\newtheorem{corollary*}[theorem*]{Corollary}
\newtheorem{proposition}[def1]{Proposition}

\theoremstyle{definition}
\newtheorem{Sdef1}{Definition}[section]

\theoremstyle{theorem}
\newtheorem{Slemma}[Sdef1]{Lemma}
\newtheorem{Scorollary}[Sdef1]{Corollary}

\newtheorem{Stheorem*}{Theorem}

\DeclareMathOperator{\Hom}{Hom}
\DeclareMathOperator{\HH}{HH}
\DeclareMathOperator{\HHcont}{HH_{cont}}
\DeclareMathOperator{\HC}{HC_{cont}}
\DeclareMathOperator{\HCsum}{HC_{cont}^{\oplus}}
\DeclareMathOperator{\HCprod}{HC_{cont}^{\Pi}}
\DeclareMathOperator{\Dctr}{D^{ctr}}
\DeclareMathOperator{\Dco}{D^{co}}
\DeclareMathOperator{\Dctrpc}{D^{ctr}_{pc}}

\newcommand{\D}[2]{ \mathrm{D}^{\mathrm{#1}}_{\mathrm{#2} }}
\newcommand{\HCa}[2]{ \mathrm{HC}^{\mathrm{#1}}_{\mathrm{#2} }}
\newcommand{\HCab}[2]{ \mathrm{HC}^{#1}_{\mathrm{#2} }}
\newcommand{\HHa}[2]{ \mathrm{HH}^{\mathrm{#1}}_{\mathrm{#2} }}
\newcommand{\Exta}[2]{ \mathrm{Ext}^{\mathrm{#1}}_{\mathrm{#2} }}

\begin{document}
\newpage
\maketitle
\begin{abstract}
	We define the appropriate homological setting to study deformation theory of complete locally convex (curved) dg-algebras based on Positselski's contraderived categories. We define the corresponding Hochschild complex controlling deformations and prove formality theorems for the Fréchet algebras of smooth functions on a manifold, the de Rham algebra and for the Dolbeault algebra of a complex manifold. In the latter case, the continuous Hochschild DGLA is $L_\infty$-equivalent to Kontsevich's extended deformation complex, to Gualtieri's deformation complex of $X$ viewed as generalized complex manifold and the (discrete) Hochschild DGLA of the derived category in case $X$ is a smooth projective variety.  We also compute the continuous Hochschild cohomology for various categories of matrix factorisations. Continuous Hochschild cohomology is compared with Hochschild cohomology of the second kind in the sense of Polishchuk-Positselski.
\end{abstract}

\tableofcontents



\section{Introduction}
For the algebra of smooth functions, Kontsevich \cite{Kontsevich} defined a dg-Lie Hochschild cocomplex $\text{HC}_{\text{PD}}(C^\infty(M))$ consisting of polydifferential operators. His famous formality result establishes an $L_\infty$-quasi isomorphism between this complex and its cohomology Lie complex which is the space of polyvectorfields on $M$ endowed with the Gerstenhaber bracket. As a corollary, one obtains that any Poisson manifold admits a (formal) deformation quantization. Similar results have been obtained for (sheaves of) polydifferential operators on a complex manifold $X$ in \cite{Damien}.
This theorem is very powerful but it has its limitations:
\begin{enumerate}
	\item It is not clear what the homological meaning of $\text{HC}_{\text{PD}}(C^\infty(M))$ is. 
	\item Every cochain $\mu$ in $\text{HC}_{\text{PD}}(C^\infty(M))$ is necessarily local. More precisely, the value of $\mu(f_1,...,f_n)$ at a point $p \in M$ only depends on the jets of $f_1,...,f_n$ in $p$. In the context of deformation quantization however, it is not possible to integrate out (non-trivial) formal deformations in a local manner. Thus, we need to enlarge $\text{HC}_{\text{PD}}(C^\infty(M))$.
\end{enumerate}
We propose to replace the Hochschild cocomplex of polydifferential operators by the Hochschild cocomplex of continuous cochains. In that framework, both of the above problems are solved. Pflaum \cite{Pflaum} has shown that there is a quasi-isomorphism of complexes
\[ \Gamma(\Lambda^* TM)[-1] \to \HC(C^{\infty}(M))\]
using rather involved geometric globalisation techniques. The right side here denotes the space of continuous Hochschild cochains for the Fréchet algebra $C^{\infty}(M)$. Combining this result with Kontsevich's formality theorem for the polydifferential Hochschild complex, it is easy to see that the above can be extended to an $L_\infty$-quasi-isomorphism; however, we did not find a reference for this consequence in the literature.
We will provide the correct homological framework to yield an interpretation of the right hand side, give a simpler, sheaf-theoretic proof of this result and in particular show how it can be extended to an $L_\infty$-quasi-isomorphism using the results of \cite{Kontsevich}. 

Let us now explain the main motivation for this paper.
For a smooth projective variety $X$, it is well-known that 
\[ \text{Ext}^2_{\mathcal O_{X \times X}}(\Delta_* \mathcal O_X, \Delta_* \mathcal O_X) = H^0(X, \Lambda^2 TX) \oplus H^\text{1}(X, TX) \oplus H^\text{2}(X, \mathcal O_X)\]
and furthermore
\[ \HH(\D{b}{coh}(X)) \cong \text{Ext}_{\mathcal O_{X \times X}}(\Delta_* \mathcal O_X, \Delta_* \mathcal O_X).\]
The group $H^1(X, TX)$ parametrizes (first order) deformations of $X$ as a scheme and, due to Kodaira \cite{Kodaira}, also its deformations as a complex manifold. Thus, if we try to understand the remaining piece $H^0(X, \Lambda^2 TX) \oplus H^\text{2}(X, \mathcal O_X)$ geometrically, we necessarily need to leave the world of complex manifolds (or commutative schemes). One solution of this problem is within the framework of Generalized Geometry \cite{gualtieri}. Gualtieri developed its deformation theory, and showed that in the case $X$ is a complex manifold (which is a special kind of generalized complex manifold), the deformations of $X$ as a generalized complex manifold are controlled by the full space 
\[H^0(X, \Lambda^2 TX) \oplus H^\text{1}(X, TX) \oplus H^\text{2}(X, \mathcal O_X).\]
He also shows that there is an analytic obstruction map whose kernel parametrizes the first order deformations that may be integrated to geometric deformations. In this paper, we will develop the algebraic counterpart. We will see that $H^0(X, \Lambda^2 TX) \oplus H^1(X, TX) \oplus H^2(X, \mathcal O_X)$ is isomorphic to the second continuous Hochschild cohomology of the Dolbeault-algebra of $X$.
More precisely, our main result reads as follows.
\begin{theorem**}
	\label{CtrDolbeault}
	Let $X$ be a complex manifold and denote by $\mathcal A(X)$ the Dolbeault algebra of $X$ viewed as a Fréchet algebra.
	The HKR-map 
	\[ \Lambda T^{1,0} X \otimes \mathcal A(X)[-1] \to \HC(\mathcal A(X))\]
	is a quasi-isomorphism. Moreover, this map can be extended to an $L_\infty$-quasi-isomorphism between $\Lambda T^{1,0} X \otimes \mathcal A(X)[-1]$ and $\HC(\mathcal A(X))$.
\end{theorem**}
This result in particular implies that the formal deformation theory of $X$ as a generalized complex manifold is equivalent to the formal continuous deformation theory of the Dolbeault algebra of $X$.
The \emph{continuous Hochschild complex} appearing in the above result is naturally understood in the context of contraderived categories of complete locally convex dg-algebras building on the work of \cite{PosiExact}. In particular, the Hochschild cocomplex then has a natural interpretation as \emph{continuous contraderived} endomorphisms of $A$ as an $A^{op} \hat \otimes A$-module. A more complete description of the main results will be given in the outline \ref{outline}.

\subsection{Outline and statement of the results}
\label{outline}
In section 1 we develop the appropriate homological framework to deal with derived categories in a continuous context. All modules and algebras considered here are modelled on complete locally convex topological vector spaces. We will work in the so-called 'split'-exact setting. A sequence of $A$-modules
\[ 0 \to M \to N \to T \to 0\]
is called \emph{exact} if it admits continuous $\mathbb K$-linear splittings; we are thus naturally in the setting of relative homological algebra. $\mathbb K$ will always either denote the field of real or complex numbers. Inspired by the work of \cite{PosiExact}, we define the contraderived category of a complete locally convex dg-algebra $A$ as the Verdier quotient
\[ \Dctr(A) := H^0(A\text{-Mod}) / \text{Ac}^{\text{ctr}}(A).\]
where $\text{Ac}^{\text{ctr}}(A)$ contains totalizations of short exact sequences and is closed under the operations of taking products, shifts and cones, see Definition \ref{ContraderivedDef}. In contrast to the discrete case, thanks to special properties of the projective tensor product \ref{ProjTensor}, we obtain a canonical dg-model for this category:
\begin{theorem*}[\ref{TheoremEnhancement}]
	Let $A$ be a clct dg-algebra. We denote by $A\text{-Mod}_{gproj}$ the full dg-subcategory of $A\text{-Mod}$ consisting of graded-projective objects. Then the functor 
	\[ H^0(A\text{-Mod}_{gproj}) \to \Dctr(A)\]
	is a triangulated equivalence.
\end{theorem*}
In subsection 1.4 we discuss why contraderived categories behave well with respect to descent. Specifically, we prove
\begin{proposition*}[\ref{LocalReductionProp}]
	Let $\mathcal A$ be a sheaf of clct dg-algebras. Let $M$ be a clct dg-module over $\mathcal A(X)$ and consider the sheaf $\mathcal M(U) := M \hat \otimes^{L}_{\mathcal A(X)} \mathcal A(U)$. Suppose that $\mathcal A$ satisfies strong descent. Then $M$ is contraacyclic if and only if there exists a cover of finite dimension $(U_i)_{i \in I}$ such that each $\mathcal M(U_i)$ is contraacyclic as an $\mathcal A(U_i)$-module (or equivalently as an $\mathcal A(X)$-module). \\ \\
	In particular, if the underlying topological space has finite covering dimension then $M$ is contraacyclic if and only if for each $x \in X$ there exists a neighbourhood $U$ of $x$ such that $\mathcal M(U)$ is contraacyclic.
\end{proposition*}
A particular case of a sheaf of clct dg-algebras $A$ satisfying the assumption of strong descent is when its $0$-th component $A^0$ is fine. \\The following subsection 1.5 introduces the main object we are interested in: The continuous Hochschild cohomology of a complete locally convex dg-algebra $A$:
\[\HH_{\text{cont}}^n(A,A) := \Hom_{\Dctr(A^e)}(A[n], A).\]
Moreover, in contrast to the discrete contraderived setting, there is a canonical $B_\infty$-Hochschild cocomplex calculating this cohomology.
\begin{theorem*}[\ref{IsDgLie}]
	Let $A$ be a clct dg-algebra. Then the continuous Hochschild complex $\HC(A, A)$ naturally carries the structure of a $B_\infty$-algebra and its cohomology is equal to $\HHa{\text{$n$}}{cont}(A,A)$.
\end{theorem*}
It will be explained in the follow-up paper that the underlying dg-Lie complex of $\HC(A, A)$ controls the deformations of $A$ as a clct curved $A_\infty$-algebra.

In \ref{CoalgebraRelations} we show how our constructions can be viewed as generalizing the theory of coderived categories of (discrete) dg-coalgebras to the setting of complete locally convex dg-algebras. To see this, recall that the category of pseudo-compact vector spaces is anti-equivalent to that of discrete vector spaces. Pseudo-compact spaces are particular cases of complete locally convex spaces; by definition a pseudo-compact space is an inverse limit of finite-dimensional spaces in the category of complete locally convex spaces. This anti-equivalence exchanges the discrete tensor product with the projective tensor product. A (discrete) dg-coalgebra is thus equivalently a pseudo-compact dg-algebra, and in fact for a dg-coalgebra $C$, there is an anti-equivalence between $\Dco(C)$ (see \cite{posi}) and the subcategory of the contraderived category of the pseudo-compact dual $C^*$ consisting of pseudo-compact modules $\Dctrpc(C^*)$. In particular, the notion of coHochschild cohomology and homology for a discrete dg-coalgebra may be interpreted in our terms as Hochschild cohomology and homology for the dual pseudo-compact dg-algebras. coHochschild complexes are important for studying algebraic models for string operations and free loop spaces, see \cite{RiveraString} and \cite{Rivera}. This example also shows that the Hochschild homology and cohomology of a clct dg-algebra is not an invariant under quasi-isomorphisms, and not even under so called \emph{Morita II-equivalences} (see \cite{HC2}) which is explained in Remark \ref{RemarkSimplicial}.

In the last subsection, we compare our construction of the contraderived category with Positselski's discrete contraderived category. In particular, we will see that there is a natural comparison functor. In the case of discrete algebras with countable bases, we show:
\begin{corollary*}[\ref{CorComparison}]
	Let $A$ be a clct cdg-algebra endowed with the finest (complete locally convex) topology and suppose that $A$ has a countable basis. Then 
	\[ \HHa{PP}{}(A,A) \cong \Hom_{\Dctr(A^e)}(A,A) =: \HHcont(A,A)\]
	where the left-hand side denotes Hochschild cohomology in the sense of Polishchuk-Positselski \cite{Polishchuk} and 
	$\Dctr(A^e)$ denotes the continuous contraderived category of $A^e$.
\end{corollary*}
We investigate the Hochschild complexes in 5 cases:
\begin{enumerate}
	\item The Fréchet algebra $C^{\infty}(M)$ of smooth functions on a real, smooth manifold $M$.
	\item The Fréchet Dolbeault dg-algebra $\mathcal A(X)$ of $(0,q)$-forms on a complex manifold $X$, together with a  '$B$-Field' $B^{0,2} \in \mathcal A^{0,2}(X)$.
	\item The Fréchet de-Rham dg-algebra $\Omega(M)$ on a real, smooth manifold $M$.
	\item For Matrix factorisations in the case of the algebras $\mathcal A^{0,*}(X)$, the discrete polynomial algebra in $n$-variables, the pseudo-compact algebra of formal power series in $n$-variables and the Fréchet-algebra of smooth functions on a smooth manifold.
	\item The pseudo-compact algebra of singular cochains $C^*(X)$ on a simplicial set $X$.
\end{enumerate}
In the spirit of \cite{Kontsevich} and \cite{Damien}, we will prove formality theorems for the Hochschild dg-Lie complexes in the first three cases, replacing geometric globalisation arguments as used in \cite{Pflaum} with elementary, sheaf-theoretic techniques. In the case of matrix factorisations, we just compute the continuous Hochschild cohomology. These calculations prove that in our framework, computing Hochschild cohomology also of curved dg-algebras is simple and straightfoward, e.g. in the case of a complex manifold, we obtain: 
\begin{theorem*}[\ref{CtrDolbeault}]
	Let $X$ be a complex manifold.
	The HKR-map 
	\[ \Lambda T^{1,0} X \otimes \mathcal A(X)[-1] \to \HC(\mathcal A(X))\]
	is a quasi-isomorphism. Moreover, this map can be extended to an $L_\infty$-quasi-isomorphism between $\Lambda T^{1,0} X \otimes \mathcal A(X)[-1]$ and $\HC(\mathcal A(X))$.
\end{theorem*}
We end the paper by giving an outlook towards Hochschild complexes of continuous enhanced dg-categories. In particular, there is a $B_\infty$-invariance result for the continuous Hochschild complex which implies that the continuous Hochschild complex of a cdg-algebra agrees with that of its continuous dg-category of finitely generated graded-projective modules. As an example, we state:
\begin{theorem*}[\ref{FormalityEnhancement}]
	Let $X$ be a compact complex manifold. Then the bounded derived category $\D{b}{coh}(X)$ has a natural clct dg-enhancement $\D{ctr,b}{}(\mathcal A(X))$ where $\mathcal A(X)$ is the Dolbeault algebra. Moreover, there is an $L_\infty$-quasi-isomorphism
	\[ \HC(\D{ctr,b}{}(\mathcal A(X)) \simeq \Lambda T^{1,0} X \otimes \mathcal A(X)[-1].\]
	The right-hand side is isomorphic to the deformation complex of $X$ as a generalized complex manifold.
	If $X$ is not compact, then the above statement holds if we replace $\D{b}{coh}(X)$ by the category of globally bounded perfect complexes on $X$.
\end{theorem*}
Similarly, we can now compute the (continuous) Hochschild cohomology of the categories of matrix factorisations for the algebras mentioned above. These results will be stated and proven carefully in the follow-up paper.

\subsection{Relations to other work}
In recent years, multiple variants of Hochschild (co-)homology have been introduced. In the case of coalgebras, there exists coHochschild theories for so-called \emph{categorical coalgebras} introduced in \cite{Rivera}. The coHochschild complex of singular chains on a simplicial set then yield algebraic models for the free loop space, see \cite{Rivera} and also for defining string operations algebraically in \cite{RiveraString}. In the case of a "single-object" categorical coalgebra, their definition agrees with our definition for the dual pseudo-compact dg-algebra. For example, if $X$ is a reduced simplicial set (i.e., which only has one $0$-cell) for which all $1$-simplices are homotopically invertible then 
\[ \HHa{cont}{*}(C^*(X)) = \mathrm{coHH}_*(C_*(X)) \simeq C_*(L |X|) \]
where $C^*(X)$ denotes the pseudo-compact algebra of singular cochains dual to the discrete coalgebra $C_*(X)$ by \cite[Theorem 20]{Rivera}, see Proposition \ref{PropRivera}. In the general "many-object" case, there is a similar result yielding the notion of \emph{pseudocompact categories} which will be discussed in future work. 

The most closely related Hochschild cohomology theory is that of Polishchuk-Positselski. It agrees with our definition in the case of a discrete algebra whose underlying vector space has at most countable dimension. Viewing the classical Hochschild complex as the product totalization of a double complex, the Polishchuk-Positselski Hochschild complex is defined as its sum totalization. In certain cases, this complex computes the contraderived endomorphisms of the algebra as a bimodule over itself \cite[Section 3]{Polishchuk}. One could define a similar totalization in the continuous setting and in fact one has inclusions 
\[ \HCsum(A, A) \hookrightarrow \HC(A, A) \hookrightarrow \HCprod(A, A)\]
where we totalized using products on the right and using sums on the left side. The left inclusion is an isomorphism if $A$ is Banach, or more generally, if $A$ admits a seminorm which is Hausdorff. For example if $A$ is the Dolbeault algebra of a compact complex manifold, then the sup-norm defines such a Hausdorff seminorm (even though it does not generate its topology). Moreover, if $A$ is concentrated in degree $0$, then all complexes are isomorphic. 

Another related Hochschild cohomology theory is called \emph{Hochschild cohomology of the second kind} following \cite{HC2} which is defined as derived endomorphisms of $A$ in the \emph{compactly generated derived category of the second kind} which was first introduced in \cite{KDM2} building on the earlier work \cite{RHHolstein}. In certain cases, this cohomology controls the (classical) deformations of the dg-category of finitely generated two-sided twisted complexes of $A$, see \cite[Corollary 7.3]{HC2}. The main example is the Dolbeault algebra of a complex algebraic manifold. In this case, the compact objects in the compactly generated derived category of the second kind is equivalent to the derived category $\D{b}{\text{coh}}(X)$ and its Hochschild cohomology of the second kind equals $\text{Ext}^*_{\mathcal O_{X \times X}}(\Delta_* \mathcal O_X, \Delta_* \mathcal O_X)$ by \cite[Theorem 7.7]{HC2}.

The homological algebra used in this paper, in particular the notion of contraderived categories based on complete locally convex spaces, is typically not used in the context of complex geometry. However, this viewpoint solves the so-called 'curvature-problem'. Consider for example a smooth projective variety $X$. By a result by \cite{Bondal}, the derived category of $X$ admits a compact generator $G$ and hence $\D{b}{\text{coh}}(X) \cong \text{Perf}(A)$ for the $A_\infty$-algebra $A = \text{RHom}(G,G)$. Now by the results of Keller \cite{Keller}, the deformation theory of $\D{b}{\text{coh}}(X)$ is equivalent to that of $A$. However, the Hochschild complex of $A$ controls \emph{curved} $A_\infty$-deformations of $A$. If the resulting deformed algebra $A'$ has non-trivial curvature, then the derived category of $A'$ is not defined and so classically, $A'$ is not a 'geometric' object. Curved (first oder) deformations of $X$ correspond to $H^2(X, \mathcal O_X)$ .which parametrises $\mathcal O_X$-gerbes and it is clear how to associate to such a class a derived category: View $[H] \in H^2(X, \mathcal O_X)$ as a \u Chech 2-cocycle which associates, with respect to a given cover, to a triple intersection $U_i \cap U_j \cap U_k$ a section $H_{ijk} \in \mathcal O(U_i \cap U_j \cap U_k)$. Then $\exp(H_{ijk}) \in \mathcal O(U_i \cap U_j \cap U_k)^\times$ twists the cocycle condition for modules over $\mathcal O_X$. In this way we obtain the derived category of sheaves over $\mathcal O_X$ twisted by the gerbe $\exp{(H_{ijk})}$. Thus, curved deformations do give rise to geometric objects, but are not realised via classical homological algebra. It is however realised through derived deformation theory for contraderived categories: The class $[H] \in H^2(X, \mathcal O_X)$ is equivalently an element in the second cohomology $H^2(\mathcal A^*(X))$ of the Dolbeault algebra of $X$, which as we show in this paper, is part of the continuous contraderived deformation complex of $\mathcal A^*(X)$. Such a cocycle $B \in \mathcal A^2(X)$ gives rise to the \emph{curved} algebra
\[ (\mathcal A^*(X), \wedge, \bar \partial, B).\]
The contraderived category, in contrast to the classical derived category, is defined for this algebra. Moreover, by \cite{antweiler}, its full subcategory given by finitely generated graded-projective modules, is equivalent to the derived category of twisted sheaves. Thus, the continuous contraderived deformation theory of $A(X)$ gives rise to interesting geometric objects, as opposed to classical deformation theory. Let us also mention that Dolbeault models can give rise to non-formal non-commutative complex manifolds, see for example \cite{PolishchukSchwarz} and $\cite{Block}$. Furthermore, T-duality for Tori gives rise to equivalences between curved Tori on the one hand, and dual non-commutative Tori on the other \cite{Block22}. To the authors knowledge, no such construction has yet been achieved in the condensed or bornological setup, which treat the sheaf $\mathcal O_X$ as the fundamental object and do not account for non-local phenomena. 

\subsection{Acknowledgements}
The author has benefitted from many discussions with his supervisor, Julian Holstein. This work has been supported and funded by the Deutsche
Forschungsgemeinschaft (DFG, German Research Foundation) – SFB 
1624 “Higher structures, moduli spaces and
integrability” – Projektnummer 506632645.

\section{Continuous contraderived categories}
\subsection{Notations and conventions}
Most algebraic objects considered in this paper are based on complete locally convex vector spaces. In essence, a complete locally convex vector space is a complete topological vector space whose topology is generated by a family of seminorms. A reader unfamiliar with the subject can for example consult \cite{Jarchow} for background although most facts we need are summarized in the Appendix. We will use the abbreviation \emph{clctvs} for a complete locally convex vector space.

We will frequently make use of the \emph{projective completed tensor product} $\hat \otimes$ (Definition \ref{projTensorProduct}) defined for complete locally convex vector spaces,  One should be aware that the category of complete locally convex spaces is \emph{not} self-enriched and the projective tensor product $\hat \otimes$ does \emph{not} preserve arbitrary coproducts. 

If $f : M \to N$ is a continuous linear map between complete locally convex spaces we denote by $\text{ker}(f)$ the kernel of $f$ taken in the category of complete locally convex spaces and by $\text{coker}(f)$ the cokernel in the same category. 
 
Let $M$ be a complete locally convex vector space and a seminorm $||\cdot||_i$ on $M$, we denote by $M_i$ the Banach space which is the completion of $M$ with respect to the seminorm $||\cdot||_i$, see Definition \ref{definitionLocalisation}.

Given a dg-algebra $A$, we denote by $A^\#$ its underlying graded algebra (where we forget the differential). Similarly, any dg-module $M$ over $A$ gives rise to a graded module $M^\#$ over the graded algebra $A^\#$. 

We will describe triangulated derived categories in this paper as homotopy categories (or quotients thereof) of (pre-)triangulated \emph{dg-categories}. A \emph{dg-category} is by definition a category enriched over (discrete) cochain-complexes, see \cite{KellerDG}. If $\mathcal C$ is a dg-category, we will denote by $\text{Hom}^*(M,N)$ the cochain-complex of homomorphisms from $M$ to $N$ and by $\text{Hom}(M,N) := H^0(\text{Hom}^*(M,N))$ the space of closed degree $0$-morphisms up to homotopy.
 
\subsection{Homological algebra of the second kind for continuous algebras}

\begin{def1}
	\label{DefAlgebras}
	\begin{enumerate}
		\item Given a cochain complex $M^*$ of complete locally convex vector spaces, we will often view $M^*$ as the "graded" space $\prod_{n \in \mathbb Z} M^n$ with its induced seminorms. 
		\item Let $M$ and $N$ be cochain complexes of clctvs. Their \emph{(projective) tensor product} is defined to be 
		\[ (M \hat \otimes N)^n := \prod_{a+b = n} M^a \hat \otimes N^b \]
		endowed with the usual differential following Koszul conventions. For the definition of the tensor product on the right-hand side, see Definition  \ref{projTensorProduct}.
		\item A \emph{clct dg-algebra $A$} is a cochain complex of clctvs together with a multiplication map of degree $0$
		$A \hat \otimes A \to A$ subject to the usual relations (note that this involves a direct product totalization so this definition is non-standard if $A$ is unbounded).
		\item Let $A$ be a clct dg-algebra. A \emph{dg-left module over $A$} is a cochain complex of clctvs $M^*$ together with a (continuous) action map 
		\[ A \hat \otimes M \to M\]
		fulfilling the usual conditions. Similarly we can define a notion of right modules. The \emph{dg-category of right $A$-modules} will be denoted by $A\text{-Mod}$. The cochain complex of morphisms $\text{Hom}_{A}^*(M,N)$ is defined to be the subspace of all homogeneous maps $\text{Hom}_{\mathbb K}(M,N)$ which commute with the $A$-action. The differential on $\text{Hom}_{A}^*(M,N)$ is given by 
		\[ d(f) = d_N \circ f - (-1)^{|f|} f \circ d_M.\]
		\item 
		Let $M$ be a right clct dg-module and $N$ a left clct dg-module over a clct dg-algebra $A$. We define their \emph{(completed/projective) tensor product} as the cokernel 
		\begin{nalign}
			\label{defTensorProduct}
			M \hat \otimes_A N := \text{coker}(M \hat \otimes A \hat \otimes N \to M \hat \otimes N)
		\end{nalign}
		where the map is given by the difference of the two actions, as usually.
	\end{enumerate}
\end{def1}
The starting point of our theory is the fact that products commute with tensor products by Lemma \ref{ProjTensor} which is of course false in non-topological settings.
\begin{corollary}
	\label{CorTensor}
	Let $M_{k \in K}$ be a family of right clct dg-modules and $N$ a left clct dg-module.
	Then there is a canonical isomorphism
	\[ (\prod_{k \in K} M_k) \hat \otimes_A N = \prod_{k \in K} (M_k \hat \otimes_A N).\]
	\begin{proof}
		This follows from Lemma \ref{ProjTensor} and that products commute with cokernels in the category of complete locally convex vector spaces by Lemma \ref{ProdCokernels}.
	\end{proof}
\end{corollary}
\begin{def1}
	Let $A$ be a clct dg-algebra. A clct dg-module $M$ is called \emph{graded-free} if its underyling graded module $M^\#$ over $A^\#$ is of the form
	\[ M^\# = F \hat \otimes_{\mathbb C} A^\#\]
	where $F$ is a graded clct vector space.
	The category of clct vector spaces admits an exact structure where a sequence 
	\[ 0 \to A \to B \to C \to 0\]
	is called exact if it admits a continuous $\mathbb C$-linear contracting homotopy (in particular, it is exact in the usual sense). This is often called the "split"-exact structure. A map from $B \to C$ fitting into a sequence as above is called an admissible epimorphism. A map from $A \to B$ fitting into a sequence as above is called an admissible monomorphism.
	
	Now let $A$ be a clct graded algebra (without differential) then we define an exact category structure on the category of graded clct $A$-modules $A\text{-Mod}_{gr}$ (which is a linear category, all morphisms are of degree $0$). A sequence
	\[ 0 \to M \to N \to L \to 0\] is called exact if 
	if it is split-exact in each degree in the above sense (i.e, the sequence admits a $\mathbb C$-linear continuous contracting homotopy).
	
	In the sense explained in \cite[Section 4.3] {PosiExact} we obtain for any clct dg-algebra $A$ an exact dg-structure on the dg-category of right clct dg-modules $A\text{-Mod}$ over $A$. This in turn defines an ordinary exact structure on the linear category $Z^0(A\text{-Mod})$, by declaring a sequence to be exact if and only if its image in $A^\#\text{-Mod}_{gr}$ is exact.
	A clct dg-module $P$ is called \emph{graded-projective} if its underyling graded module $P^\#$ is projective with respect to the above exact-structure. In detail, this means the following. Given two graded clct modules $M$ and $N$ over $A^\#$ and an $A^\#$-linear map $M \to N$, which is an admissible epimorphism in each degree (as $\mathbb C$-linear continuous maps of clctvs) and an $A^\#$-linear map $P^\# \to N$, there is a lift to an $A^\#$-linear map $P^\# \to M$.
\end{def1}
\subsection{The contraderived category}
\begin{def1}
	\label{ContraderivedDef}
	Let $A$ be a clct dg-algebra. The exact dg-category $A\text{-Mod}$ has shifts and cones. Following \cite[Section 5.1]{PosiExact}, we define the category $\text{Ac}^{\text{abs}}(A)$ as the minimal thick full subcategory of $H^0(A\text{-Mod})$ which contains all totalizations of short exact sequences in $Z^0(A\text{-Mod})$ (w.r.t to the specified exact structure). The absolute derived category of $A$ is defined as the Verdier quotient
	\[ \text{D}^{\text{abs}}(A) := H^0(A\text{-Mod}) / \text{Ac}^{\text{abs}}(A).\]
	Infinite products exist in the dg-category $A\text{-Mod}$ and products of exact sequences are exact (one can take the product of the corresponding continuous homotopies). Then again due to \cite[Section 5.1]{PosiExact} we may define the minimal thick triangulated subcategory $\text{Ac}^{\text{ctr}}(A)$ of $H^0(A\text{-Mod})$ which contains $\text{Ac}^{\text{abs}}(A)$ and is closed under direct products and consider the Verdier quotient 
	\[ \Dctr(A) := H^0(A\text{-Mod}) / \text{Ac}^{\text{ctr}}(A).\]
	Modules in $\text{Ac}^{\text{ctr}}(A)$ will be termed \emph{contraacyclic}.
\end{def1}
\begin{remark}
	In the definition of the contraderived category we demanded that \emph{arbitrary} products of contraacyclic objects are again contraacyclic. However, we would obtain the same category if we only required that \emph{countable} products of contraacyclic objects are again contraacyclic.  This is true since we only need countable products leading up to the proof of Theorem \ref{TheoremEnhancement} 

\end{remark}
\begin{remark}
	In general, Verdier quotients give rise to non locally small categories. In our case, $\Dctr(A)$ is locally small due to Theorem \ref{TheoremEnhancement}. It is not clear to us whether $D^{abs}(A)$ is locally small. In the case of discrete algebras, one can also define the so-called coderived category, \cite{posi}. In our setting, we cannot define it since the category $A\text{-Mod}$ does, in general, not have sums because the projective tensor product may not commute with sums.
\end{remark}
\begin{lemma}
	\label{LemmaHom}
	Let $A$ be a clct dg-algebra and $M$ a clct dg-module over $A$ which is graded-projective and $N$ be a clct dg-module.
	Then \[H^0(\Hom^*_{A\text{-Mod}}(M, N)) = \Hom_{\Dctr(A)}(M, N)\]
	\begin{proof}
		See \cite[Theorem 5.5]{PosiExact}.
	\end{proof}
\end{lemma}
The previous Lemma tells us that we can compute morphism spaces in the contraderived category by finding for any module $M$ a graded projective replacement $P \to M$. In order to verify that some map is an isomorphism in $\Dctr(A)$, here is a criterion:
\begin{lemma}
	\label{LemmaTot}
	Let $A$ be a clct dg-algebra. Suppose \[M^* = ... \to M^{-2} \to M^{-1} \to M^0 \to 0\] is a cochain complex of clct dg $A$-modules bounded from above which is exact (horizontally). Then its product totalization $\text{Tot}^{\Pi}(A)$ is contraacyclic.
	\begin{proof}
		See \cite[Theorem 5.11]{PosiExact}.
	\end{proof}
\end{lemma}

As in the discrete case, it is simple to show that:
\begin{lemma}
	\label{LemmaEquivalent}
	For a clct dg-module the following are equivalent:
	\begin{enumerate}
		\item $M^\#$ is graded-projective.
		\item $M^\#$ is a direct summand of a graded-free clct module.
		\item $M$ is a direct summand of a clct dg-module which is graded-free.
		\item $M$ is a direct summand of a clct dg-module which is graded-free and the complement is contraacyclic.
	\end{enumerate}
	\begin{proof}
		The equivalence of the first two statements follows by considering the action $M \hat \otimes A \to M$. As products commute with tensor products we obtain that the product totalization of the one-sided bar resolution
		\[ ... \to M \hat \otimes A \hat \otimes A \to M \hat \otimes A \to M\]
		is a graded-free replacement for $M$. Consider the pushforward
		\[ \epsilon_* : \Hom^*(M, \text{bar}(M)) \to \Hom^*(M,M)\]
		along the augmentation $\epsilon : \text{bar}(M) \to M$. By Lemma \ref{LemmaHom}, we see that $\epsilon_*$ is a quasi-isomorphism and since $\text{bar}(M) \to M$ is an admissible epimorphism and $M$ is graded-projective, we furthermore obtain that $\epsilon_*$ is surjective. Any surjective quasi-isomorphism of complexes is surjective on cocycles. We thus obtain a right-inverse of $\epsilon$ and so $1 \Rightarrow 4$. The other implications are obvious.
	\end{proof}
\end{lemma}

In contrast to the discrete case, we have:
\begin{corollary}
	\label{ProdProj}
	Let $A$ be a clct dg-algebra. Let $P_{n \in N}$ be an arbitrary family of graded-projective clct $A$-modules.
	Then their product $\prod_{n \in N} P_n$ is graded-projective.
	\begin{proof}
		First suppose that $P_n$ is graded-free. Then we can write $P_n^\# = F_n \hat \otimes A^\#$. We have seen that 
		\[ \prod_{n \in N} P_n^\# = (\prod_{n \in N} F_n) \hat \otimes A^\#\]
		which is graded-free.
		By the above Lemma we obtain the same for arbitrary graded projective $P_n$.
	\end{proof}
\end{corollary}

\begin{theorem}
	\label{TheoremEnhancement}
	Let $A$ be a clct dg-algebra. We denote by $A\text{-Mod}_{gproj}$ the full dg-subcategory of $A\text{-Mod}$ consisting of graded-projective objects. Then the functor 
	\[ H^0(A\text{-Mod}_{gproj}) \to \Dctr(A)\]
	is a triangulated equivalence.
	\begin{proof}
		This follows from \cite[Theorem 5.10]{PosiExact} using Corollary \ref{ProdProj} and by the existence of the (graded-free) bar resolution as in the proof of Lemma \ref{LemmaEquivalent}.
	\end{proof}
\end{theorem}
\begin{remark}
	We will see in the follow-up paper, that there exists a model structure on the category $Z^0(A\text{-Mod})$ whose homotopy category yields the contraderived category of $A$. The fibrations are the admissible epimorphisms and the cofibrations are admissible monomorphisms with graded-projective cokernel. 
\end{remark}
\begin{remark}
	The category $\Dctr(A)$ is very large as arbitrary complete locally convex spaces appear as fibers of graded-free modules. By Theorem \ref{TheoremEnhancement}, it is natural to consider subcategories generated by objects $(V \hat \otimes A, d)$ for $V$ an element of a specified \emph{set} of complete locally convex spaces (e.g., finite-dimensional spaces, pseudo-compact spaces, seperable Hilbert spaces, Hilbert spaces with Schauder basis of bounded cardinality,...). It would be interesting to study the $K$-theory of such subcategories. In fact, localising at all objects of the form $(V \otimes A,d)$ with $V$ finite-dimensional yields the \emph{compactly generated derived category of the second kind}, \cite{KDM2}. If we only localise at the module $A$, we obtain the usual (discrete) derived category of $A$.
\end{remark}
\begin{def1}
	\label{Enhancement}
	As seen in the last Theorem, we have a canonical dg-enhancement of $\Dctr(A)$. We will denote this dg-category by
	\[ \text{\underline D}^{\text{ctr}}(A) := A\text{-Mod}_{gproj}.\]
\end{def1}
\begin{remark}
	If $A$ is a Fréchet dg-algebra we may restrict all dg-modules to be also Fréchet and may develop a completely analogous theory (because \emph{countable} products of Fréchet spaces are Fréchet) to obtain the contraderived category $\D{ctr}{Fré}(A)$, without ever having to talk about general clctvs. Moreover, there is a functor
	\[\D{ctr}{Fré}(A) \to \Dctr(A)\] 
	which is fully faithful because each Fréchet dg-module admits a graded-projective resolution (as we use products for totalization!) which is also Fréchet. Thus, we may also describe $\D{ctr}{Fré}(A)$ as the homotopy category of graded-projective Fréchet dg-modules. This is very convenient, as it allows one to remain in the analytically nice realm of Fréchet spaces, and on the other hand it has a good homological interpretations. In particular, this in contrast to the usual setup of homological algebra (also called homological algebra of the first kind) where the necessity of direct sum totalization forces one to leave the world of Fréchet spaces, and for example work with bornological modules. 
\end{remark}
\begin{def1}
	Let $f : A \to B$ be a map of clct dg-algebras.
	\begin{enumerate}
		\item 
		The tensor product over $A$ induces an extension of scalar functors of the dg-categories of graded-projectives 
		\[(-) \hat \otimes_A B :  \text{\underline D}^{\text{ctr}}(A) \to \text{\underline D}^{\text{ctr}}(B)\]
		(see Definition \ref{Enhancement}) by Corollary \ref{CorTensor}
		and thus a derived functor 
		\[ Lf_* = (-) \hat \otimes_A^{L} B : \Dctr(A) \to \Dctr(B)\]
		which may be described by first choosing graded-projective replacements and then applying $(-) \hat \otimes_A B$ or alternatively by choosing a graded-projective replacement for $B$.
		In the case that $B$ is graded-flat over $A$ (for example if $B$ is graded-projective), then the tensorproduct is already derived so we obtain that 
		\[ f_* = (-) \hat \otimes_A B : \Dctr(A) \to \Dctr(B)\]
		is well-defined.
		\item There is also a restriction of scalars functor:
		\[ f^* : H^0(B\text{-Mod}) \to H^0(A\text{-Mod}).\]
		By definition of the exact structure, this functor preserves totalizations of short exact sequences. Clearly, it is also compatible with the triangulated structures. Moreover, it sends products to products. It follows that we also obtain a functor on the Verdier quotient:
		\[  f^* : \Dctr(B) \to \Dctr(A).\]
	\end{enumerate}
\end{def1}
\begin{lemma}
	\label{Adjunction}
	Let $f : A \to B$ be a map of clct dg-algebras. Then $Lf_*$ is left adjoint to $f^*$. 
	\begin{proof}
		Let $P$ be a graded-projective $A$-module and $M$ a $B$-module. There is a canonical isomorphism of cochain complexes:
		\[\Hom_A^*(P, f^*M) = \Hom_B^*(P \hat \otimes_A B, M). \]
		Since $P$ is graded-projective, the left side computes $\Hom_{\Dctr(A)}(P, f^*M)$.\\ $P \otimes_A B$ is graded-projective over $B$ so the right-hand side computes\\ $\Hom_{\Dctr(B)}(Lf_* P, M)$. Moreover, every module admits a graded-projective replacement by Theorem \ref{TheoremEnhancement}.
	\end{proof}
\end{lemma}
\subsubsection{Non-negatively graded clct-dg algebras}
Suppose $A$ is concentrated in non-negative degrees. Then there is a map of algebras $\pi_0 : A \to A^0$ and we obtain a corresponding extension by scalars functor $(-) \hat \otimes^{L}_A A^0 : \Dctr(A) \to \Dctr(A^0)$.
\begin{proposition}
	Suppose $A$ is concentrated in non-negative degrees. Then the extension of scalars functor $(-) \hat \otimes^{L}_A A^0 : \Dctr(A) \to \Dctr(A^0)$ reflects isomorphisms. In other words, if $N$ is graded-free $N = (W \hat \otimes A^\#, d_0+d_1+...)$, then $N$ is a contraacyclic $A$-module if and only if $(W \hat \otimes A^0, d_0)$ is a contraacyclic $A^0$-module.
	\begin{proof}
		Suppose $N$ is graded-free as above such that $(W \hat \otimes A^0, d_0)$ is a contraacyclic $A^0$-module. Let $M = (V \hat \otimes \mathcal A, d^M_0 + d^M_1 + ...)$ be any other graded-free module. The silly filtration $\sigma_{\geq k} A$ is an $A$ bimodule for each $k$ and there are $A$-bimodule maps $\sigma_{\geq k+1} A \to \sigma_{\geq k}$ (all of this assumes $A$ to be non-negatively graded).
		The complex $\Hom^*_{A}(M,N)$ is naturally decreasingly filtered as
		\[ F^k \Hom^*_{A}(M,N) =  \Hom^*_{A}(M,N \hat \otimes_{A} \sigma_{\geq k} A).\]
		The filtered pieces are then isomorphic to $\Hom^*_{A^0}((V \hat \otimes A^0, d_0),(W \hat \otimes A^k, d_0))$. As $V \hat \otimes A$ is totalized using products, this filtration is complete. Moreover it is bounded above ($ F^0 \Hom^*_{A}(M,N) = \Hom^*_{A}(M,N)$). The corresponding spectral sequence thus converges. By assumption we have that $(W \hat \otimes A^0, d_0)$ is contraacyclic. But then \[W \hat \otimes A^0 \hat \otimes_{A^{0}} A^i = W \hat \otimes A^i\] is also contraacyclic (we do not need to derive as $W \hat \otimes A^0$ is graded-free). This shows that $\Hom^*_{A^0}((V \hat \otimes A^0, d_0),(W \hat \otimes A^k, d_0))$ has vanishing cohomology and by the spectral sequence we deduce acyclicity of $ \Hom^*_{A}(M,N)$. As $M$ was arbitrary, we deduce $N$ to be contraacyclic.
	\end{proof}
\end{proposition}
If we consider the Dolbeault algebra $\mathcal A^{0,*}(X)$ of a complex manifold $X$ and view a graded-free module $M^\# = V \hat \otimes \mathcal A^\#$ as a flat partial superconnection on the trivial vector bundle with (infinite-dimensional) fibre $V$ then this proposition can be interpreted as saying $M$ is contraacyclic if and only if the underlying complex of smooth vector bundles $(V \hat \otimes C^\infty(M), d_0)$ is contraacyclic.

\subsection{Limits and colimits}
Let $A$ be a clct dg-algebra. The dg-category of $A$-modules has products but may fail to have sums because the projective tensor product does in general not commute with sums. However, the contraderived category has direct sums.
\begin{lemma}
	The category $\Dctr(A)$ has (small) products and coproducts.
	\begin{proof}
		The statement for existence of products easily follows as cohomology commutes with products and we can compute Hom's using graded projective replacements. The underlying clct-dg vector space of a product $\prod A_i$ of clct dg-modules is simply the usual direct product in the category of clct vector spaces.\\ \\
		For direct sums we need to be a bit more careful. Let $A_i$ be a family of objects in $\Dctr(A)$. We may assume that $A_i$ are graded-free, i.e, $A_i^\# = V_i \hat \otimes A^\#$. The underlying clct graded $A^\#$ module of their coproduct is defined as 
		\[ (\oplus A_i)^\# := (\oplus V_i) \hat \otimes A^\#\]
		where $\oplus V_i$ is the (categorical) direct sum operation in the category of clct dg-vector spaces. The differential on any graded-free clct-dg module $W \hat \otimes A^\#$ is determined by the restrictions
		\[ d^k : W \to W \hat \otimes A^k, k \in \mathbb Z.\]
		We denote these differentials as $d^k_i$ for each $A_i$. Then the differentials $d^k :  (\oplus V_i) \to  (\oplus V_i) \hat \otimes A^k$ on $\oplus A_i$ are equivalently specified by a family $(d^k)_i$ given by the composition 
		\[ (d^k)_i = (\iota_i \hat \otimes id) \circ d^k_i : V_i \to V_i \hat \otimes A^k \hookrightarrow (\oplus V_j) \hat \otimes A^k\]
		It is clear that the differential then squares to $0$. Then, forgetting differentials for a moment, we have for any other dg-module $M$
		\begin{nalign}
			\Hom^*_{A}((\oplus V_i) \hat \otimes A, M) = &\Hom^*_{\mathbb C}(\oplus V_i, M) \\= &\prod \Hom^*_{\mathbb C}(V_i, M) = \prod \Hom^*_{A}(V_i \hat \otimes A, M).
		\end{nalign}
		By construction of the differential on $(\oplus V_i) \hat \otimes A$ this isomorphism is clearly compatible with the differentials. Taking cohomology on both sides proves the claim.
	\end{proof}
\end{lemma}
\begin{remark}
	If $A$ is Fréchet, then $\D{ctr}{Fré}(A)$ only has countable products and finite sums.
\end{remark}
\begin{corollary}
	Let $A_1 \to A_2 \to ...$ be a diagram in $\Dctr(A)$ indexed over $\mathbb N$. Then its homotopy colimit in $\Dctr(A)$ exists. Similarly, any diagram of the form $B_1 \leftarrow B_2 \leftarrow ...$ has a homotopy limit in $\Dctr(A)$.
	\begin{proof}
		This follows as $\Dctr(A)$ is triangulated and has direct sums and direct products, see \cite{HomLimits}.
	\end{proof}
\end{corollary}
\subsection{Sheaves and contraderived categories}
\label{subsectionSheaves}
\begin{theorem}
	\label{OpenIsProjective}
	Let $M$ be a smooth manifold and $U$ an open subset of $M$. Then $C^{\infty}(U)$ is projective over $C^\infty(M)$ (with respect to the split-exact structure).
	\begin{proof}
		Due to Ogneva \cite[Theorem 2]{ogneva}, the above statement is true in case $U$ is contained in a coordinate chart. Now let $U \subseteq M$ be an arbitrary open subset. Then we can find $U_i$ which are contained in coordinate charts and which cover $U$. Choose a (locally finite) partition of unity $\epsilon_i$ subordinate to this cover. By the above mentioned result, there are continuous $C^\infty(M)$-linear sections $s_i$ to the action maps
		\[ \mu_i : C^{\infty}(M \times U_i) \to C^{\infty}(U_i).\]
		Then $\sum_{i} (1 \otimes \epsilon_i) s_i \circ r_i$ is a continuous $C^\infty(M)$-linear section of the action $\mu : C^{\infty}(M \times U) \to C^{\infty}(U)$, where $r_i$ is the restriction $r_i : C^{\infty}(U) \to C^{\infty}(U_i)$. As $C^{\infty}(M \times U)$ is free this proves the claim.
	\end{proof}
\end{theorem}
\begin{example}
	\label{ExRes}
	Let $M$ be a manifold and $U$ an open subset of $M$. Then the functor
	\[ (-) \hat \otimes_{C^{\infty}(M)} C^{\infty}(U) : \Dctr(C^{\infty}(M)) \to \Dctr(C^{\infty}(U))\]
	is already derived (by the previous Theorem).
	Consider dg-modules $M$ and $N$ over $C^{\infty}(M)$, and choose a graded-projective replacement $P \to M$. Then the presheaf $U \mapsto N \hat \otimes_{C^{\infty}(M)} C^{\infty}(U)$ is actually a sheaf by Proposition \ref{IsSheaf}. Furthermore, \\
	$U \mapsto \Hom^*_{C^{\infty}(M)}(P, N \hat \otimes_{C^{\infty}(M)} C^{\infty}(U))$ is also a complex of sheaves and its global sections compute the contraderived Hom of $M$ and $N$. Moreover, each \\$\Hom^*_{C^{\infty}(M)}(P, N \hat \otimes_{C^{\infty}(M)} C^{\infty}(U))$ is by adjunction isomorphic to \[\Hom^*_{C^{\infty}(U)}(P \hat \otimes_{C^{\infty}(M)} C^{\infty}(U) , N \hat \otimes_{C^{\infty}(M)} C^{\infty}(U)).\]
	which in turn is the contraderived Hom
	\[\Hom_{\Dctr(U)}(N \hat \otimes_{C^{\infty}(M)} C^{\infty}(U) , N \hat \otimes_{C^{\infty}(M)} C^{\infty}(U))\]
	because $C^{\infty}(U)$ is (graded-)projective over $C^{\infty}(M)$.
	This construction is useful for globalization arguments. The same constructions work for the De-Rham algebra and the Dolbeault algebra of a complex manifold.
\end{example}
\begin{def1}
	Let $X$ be a topological space and $\mathcal A$ a sheaf of clct dg-algebras. We say that 
	$\mathcal A$ satisfies strong descent if for any cover of finite dimension $(U_i)_{i \in I}$ of any open $U \subseteq X$ the natural map
	\[ \mathcal A(U) \to \check C^*((U_i)_{i \in I}, \mathcal A)\]
	is an isomorphism in the contraderived category.  
\end{def1}
\begin{lemma}
	Let $\mathcal A$ be a sheaf of clct dg-algebras such that $\mathcal A^0$ is a fine sheaf. Then $\mathcal A$ satisfies strong descent. 
	\begin{proof}
		As $\mathcal A^0$ is fine, we obtain for any finite-dimensional cover a horizontally-exact finite length complex
		\[ 0 \to \mathcal A(U) \to \check C^0((U_i)_{i \in I}, \mathcal A) \to ... \to \check C^n((U_i)_{i \in I}, \mathcal A) \to 0\]
		Here we use that a partition of unity always gives rise to continuous splittings.
		Its totalization is contraacyclic because it is bounded above (as $(U_i)$ has finite dimension).
	\end{proof}
\end{lemma}
\begin{proposition}
	\label{IsSheaf}
	Let $X$ be a topological space and $\mathcal A$ a sheaf of clct dg-algebras. Let $M$ be a clct dg-module over the global sections $\mathcal A(X)$ of $\mathcal A$. Then the presheaf
	\[ \mathcal M(U) := M \hat \otimes^{L}_{\mathcal A(X)} \mathcal A(U)\]
	is a sheaf. 
	\begin{proof}
		The product in the definition of the \v Cech-complex commutes with the tensor product.
	\end{proof}
\end{proposition}
\begin{remark}
	The above proposition is in particular interesting when each $\mathcal A(U)$ is graded-flat over $\mathcal A(M)$ (this is the case for smooth functions, the de-Rham algebra and the Dolbeault algebra) since the tensor product does not have to be derived in that case.
\end{remark}
\begin{proposition}
	\label{LocalReductionProp}
	Let $\mathcal A$ be a sheaf of clct dg-algebras. Let $M$ be a clct dg-module over $\mathcal A(X)$ and consider the sheaf $\mathcal M(U) := M \hat \otimes^{L}_{\mathcal A(X)} \mathcal A(U)$. Suppose that $\mathcal A$ satisfies strong descent. Then $M$ is contraacyclic if and only if there exists a cover of finite dimension $(U_i)_{i \in I}$ such that each $\mathcal M(U_i)$ is contraacyclic as an $\mathcal A(U_i)$-module (or equivalently as an $\mathcal A(X)$-module). \\ \\
	In particular, if the underlying topological space has finite covering dimension then $M$ is contraacyclic if and only if for each $x \in X$ there exists a neighbourhood $U$ of $x$ such that $\mathcal M(U)$ is contraacyclic.
	\begin{proof}
		
		First of all it is clear that if $M$ is contraacyclic, then all $\mathcal M(U)$ are contraacyclic as the tensorproduct is derived. Conversely if we have a finite-dimensional cover $(U_i)_{i \in I}$ such that each $\mathcal M(U_i)$ is contraacyclic it also follows that each $\mathcal M(U)$ is contraacyclic if $U \subseteq U_i$ for some $i$. We assume that $M$ is graded-projective. Then we have an isomorphism
		\[ M \to M \hat \otimes_{\mathcal A(X)} \check C^*((U_i), \mathcal A)\]
		since $\mathcal A$ satisfies strong descent. As the class of contraacyclic modules is closed under taking products we have that each $M \hat \otimes_{\mathcal A(X)} \check C^k((U_i), \mathcal A)$ is contraacyclic. As the cover has finite dimension we have the finite sequence
		\[ 0 \to M \hat \otimes_{\mathcal A(X)} \check C^0((U_i), \mathcal A) \to ... \to M \hat \otimes_{\mathcal A(X)} \check C^n((U_i), \mathcal A) \to 0\]
		consisting of contraacyclic modules. As a finite extension of contraacyclic modules, its totalization is thus also contraacyclic.
	\end{proof}
\end{proposition}

\subsection{Continuous Hochschild cohomology and homology}
Let $A$ be a clct dg-algebra. We consider the clct dg-algebra $A^e := A^{op} \hat \otimes A$ and view $A$ as a right module over $A^e$. 
\begin{def1}We define the \emph{continuous Hochschild cohomology} of $A$ with coefficients in an $A^e$-module $M$ by 
\[ \HH_{\text{cont}}^n(A,M) := \Hom_{\Dctr(A^e)}(A[n], M). \]
Similarly, the \emph{continuous Hochschild homology} of $A$ with coefficients in $M$ is defined by 
\[ \HH^{\text{cont}}_n(A,M) = H^{-n} (A \hat \otimes_{A^e}^{L} M).\]
\end{def1}
As usually, we would like to have a Hochschild cocomplex which computes Hochschild cohomology but moreover carries the structure of a dg-Lie algebra describing the (continuous) deformations of $A$. Crucially, the existence of the Lie bracket relies on continuity, and does not work for contraderived categories over discrete algebras.
\begin{def1}
	We define the dg $A^e$-module \[\text{Bar}^{-n}(A) = A^{op} \hat \otimes (A \hat \otimes ... \hat \otimes A) \hat \otimes A.\]
	We use the usual formulas induced by the action to obtain a cochain complex $\text{Bar}^{*}(A)$ of dg $A^e$-modules. Each $\text{Bar}^{-n}(A)$ is just the tensor product of $A^e$ with $A^{\hat \otimes n}$ so that each $\text{Bar}^{-n}(A)$ is graded-free and hence graded-projective. We have the same augmentation as in the classical setup and also a $\mathbb C$-linear continuous contracting homotopy ($a_0 \otimes ... \otimes a_n \mapsto 1 \otimes a_0 \otimes ... \otimes a_n$) so that 
	\[... \to \text{Bar}^{-2}(A) \to \text{Bar}^{-1}(A) \to \text{Bar}^{0}(A) \to A \to 0\]
	is a resolution of $A$ by graded-projective $A^e$-modules. Similarly, we can define the reduced version $\text{Bar}_{red}^{-n}(A) =  A^{op} \hat \otimes (A/ \mathbb K \hat \otimes ... \hat \otimes A/ \mathbb K) \hat \otimes A$. We consider the product totalization $\text{Bar}_{red}(A)$.
\end{def1}
\begin{lemma}
	Let $M$ be an $A^e$-module. Then we have
	\[ \HH_{\text{cont}}^n(A,M) = H^n(\Hom^*_{A^e}(\text{Bar}_{red}(A), M)).\]
	Similarly,
	\[ \HH^{\text{cont}}_n(A,M) = H^{-n}(\text{Bar}_{red}(A) \hat \otimes_{A^e} M).\]
	\begin{proof}
		First of all, $\text{Bar}_{red}(A) \to A$ is an isomorphism in the contraderived category by Lemma \ref{LemmaTot}.
		Moreover, $\text{Bar}_{red}(A)$ is graded-projective by Corollary \ref{ProdProj}. Thus, the result follows from Lemma \ref{LemmaHom}.
	\end{proof}
\end{lemma}
\begin{def1}
	\label{DefHochschild}
	As in the discrete case, we consider the sequence of cochain-complexes $\HC^{n,*}(A, M) := \Hom^*_{A^e}(\text{Bar}^{n+1}_{red}(A), M)$. Furthermore, we define the \emph{continuous Hochschild cocomplex of $A$ with values in $M$} as \[\HCa{*}{cont}(A, M) := \Hom_{A^e}(\text{Bar}_{red}(A), M)[-1].\]
	The \emph{continuous Hochschild complex of $A$ with values in $M$} is defined as
	\[ \HCa{cont}{*} (A,M) := \text{Bar}_{red}(A) \hat \otimes_{A^e} M.\]
	
	In the case that $A$ is concentrated in degree $0$, the definition of the Hochschild cocomplex agrees with the definition given in \cite{Pflaum}.
	Note that we may even view $\HCa{cont}{*}(A,M)$ as a complete locally convex chain complex, in contrast to the case of the continuous Hochschild cocomplex. In fact, $\HCa{cont}{*}(A,(-))$ defines a functor $\Dctr(A^e) \to \Dctr(\mathbb K)$ whereas $\HCa{*}{cont}(A, (-))$ is a functor $\Dctr(A^e)\to \text{D}(\mathbb K)$ where $\D{}{}(\mathbb K)$ is the (discrete) derived category of $\mathbb K$.
	
	In general, $\HC(A, M)$ is neither the direct sum totalization nor the direct product totalization of $\HCa{*,*}{cont}(A, M)$. We will thus also define $\HCprod(A, M)$ as the direct product totalization of $\text{HC}_\text{cont}^{*,*}(A, M)$ and  $\HCsum(A, M)$ as the direct sum totalization. We have the following inclusions:
	\[ \HCsum(A, M) \hookrightarrow \HC(A, M) \hookrightarrow \HCprod(A, M).\]
	The first map is obvious and the second map is a consequence of Corollary \ref{CorProdHom} (such a map also exists for discrete algebras, but in that case it is surjective rather than injective).
	The left inclusion is an isomorphism if $M$ is a Banach space by the following Lemma. Moreover, both maps are isomorphisms if $A$ is concentrated in degree $0$. As usually, one may define a circle product (when $M = A$):
	\begin{nalign}
		&\circledcirc : \HCab{n,p}{cont}(A, A) \otimes \HCab{m,q}{cont}(A, A) \to \HCab{n+m,p+q}{cont}(A, A)\\
		&(f \circledcirc g) = \sum_{i = 0}^n (-1)^i f \circ (\text{id}^{\otimes i} \hat \otimes g \hat \otimes \text{id}^{n-i}).
	\end{nalign}
	In the above formula, we use the interpretation that each $\HCab{n,*}{cont}(A, A)$ may be viewed as the space of $(n+1)$-multilinear, continuous maps $A^{\hat \otimes n+1} \to A$.
	As before, it is easy to show that it extends to an operation both on $\HCsum(A, A)$ and on $\HCprod(A, A)$. It also extends to an operation on $\HC(A, A)$ but this is less obvious and we will show it in the following section (and we see no reason for this to hold true for discrete algebras).
\end{def1}

\begin{theorem}
	\label{IsDgLie}
	Let $A$ be a clct dg-algebra. Then the continuous Hochschild complex $\HC(A, A)$ naturally carries the structure of a $B_\infty$-algebra. Moreover, we have inclusions of $B_\infty$-algebras
	\[ \HCsum(A, A) \hookrightarrow \HC(A, A) \hookrightarrow \HCprod(A, A).\]
	\begin{proof}
		Here we will only show that $\HC(A, A)$ is canonically a dg-Lie algebra. The full $B_ \infty$-structure can be identified similarly and will be dealt with in more general cases in the follow-up paper.
		Recall that we have the inclusions
		\[ \HCsum(A, A) \hookrightarrow \HC(A, A) \hookrightarrow \HCprod(A, A).\]
		and there are Lie brackets defined on $\HCsum(A, A)$ and $\HCprod(A, A)$ by the usual formulas. We show that the circle product on $\HCprod(A, A)$ is closed on the subspace $\HC(A, A) \subseteq \text{HC}^{\Pi}(A, A)$. We write $A$ as a reduced inverse limit of Banach spaces $A_j$ as in Definition \ref{definitionLocalisation} and denote the projections by $p_j : A \to A_j$.
		
		Thus, let $f,g \in \HC(A, A)$. By Corollary \ref{CorProdHom}, it suffices to show that \[p_j \circ (f \circledcirc g) \in \HCsum(A, A_j)\] for every semi-norm $||\cdot||_{j}$ of $A$ (technically speaking, $A_j$ is not an $A^e$-module anymore, so for all the following statements we use the interpretation as multilinear maps and forget about differentials). By definition of the circle product, we have
		\[p_j \circ (f \circledcirc g) = (p_j \circ f) \circledcirc g\]
		Now 
		\[ p_j \circ f \in \HC(A, A_j) = \HCsum(A, A_j)\]
		by Lemma \ref{LemmaProdPlus}. Thus, we may write
		\[ p_j \circ f = f_1 + ... + f_n\]
		for some $f_i \in \HCab{n_i,*}{cont}(A, A_j) = \Hom_{\mathbb C}(A^{\hat \otimes n_i+1}, A_j)$. $A_j$ is a Banach space, so by definition of the projective tensor product we may find a seminorm $s_i$ on $A$ such that $f_i$ factors through $p_{s_i} : A^{\hat \otimes n_i+1} \to A_{s_i}^{\hat \otimes n_i+1}$ as $f_i = \tilde f_i \circ p_{s_i}$. Clearly, we have $\tilde f_i \circledcirc g = f_i \circledcirc (p_{s_i} \circ g)$ for each $i$. But now the same reasoning applies to $p_{s_i} \circ g$ by which 
		\[ p_{s_i} \circ g \in \HCsum(A, A_{s_i}).\]
		and thus 
		\[ f_i \circledcirc g = \tilde f_i \circledcirc (p_{s_i} \circ g) \in \HCsum(A, A_{j})\]
		so that 
		\[ p_j \circ (f \circledcirc g) \in \HCsum(A, A_{j})\]
		for every seminorm $||\cdot||_{j}$ of $A$. As explained above, this implies \[f \circledcirc g \in \HC(A, A) \subseteq \HCprod(A, A).\]
	\end{proof}
\end{theorem}

\subsection{Pseudocompact algebras and coalgebras}
\label{CoalgebraRelations}
It is well-known that the dual of a dg-coalgebra gives rise to a dg-algebra. The converse is not true. This can be explained by the fact that every dg-colagebra is the union of its finite-dimensional subcoalgebras yet there is no analogous statement for dg-algebras. This issue may be fixed when incorporating a topology on the dual of a dg-coalgebra. Namely, the dual $A^*$ of a dg-coalgebra $A$ may be regarded as the inverse limit of all $C^*$ for $C$ a finite-dimensional subcoalgebra of $C$. A dg-algebra having this property is called a \emph{pseudo-compact dg-algebra}. It turns out that pseudo-compact dg-algebras are precisely anti-equivalent to dg-coalgebras. An introduction is available at \cite{Guan}. Note that a pseudo-compact dg-algebra is also a complete locally convex dg-algebra in a natural way (each norm on the finite-dimensional quotients gives rise to a semi-norm).
\begin{def1}
	The \emph{coderived category} of a dg-coalgebra $C$ is defined as the Verdier quotient
	\[ \Dco(C) := H^0(C-\text{Comod}) / \text{Ac}^{\text{co}}(C).\]
	where $\text{Ac}^{\text{co}}(C)$ consists of so-called coacyclic comodules. This is the minimal thick triangulated subcategory of $H^0(C-\text{Comod})$ containing totalizations of short exact sequences and closed under taking direct \textbf{sums}, see \cite{posi}.
\end{def1}
Taking duals induces an anti-equivalence $\Dco(C) \cong \Dctrpc(C^*)^{op}$. The right-hand side here denotes the opposite category of the full subcategory of the continuous contraderived category of $C^*$ consisting of all modules whose underlying topologies are pseudo-compact. Note that it precisely interchanges graded-projective $C^*$-modules and graded-injective $C$-comodules.
\begin{remark}
	$\Dctrpc(A)$ for pseudocompact $A$ can either be defined as the full subcategory of $\Dctr(A)$ consisting of pseudo-compact modules or it can be defined as the Verdier quotient of pseudocompact modules by contraacyclic pseudo-compact modules (without ever referring to general complete locally convex spaces). This is again due to the fact that the standard graded-projective resolution applied to a pseudocompact module is again pseudocompact.
\end{remark}
\begin{remark}
	Our approach can thus be understood as a generalization of coderived categories over coalgebras (or equivalently contraderived categories of pseudocompact algebras) to that of complete locally convex algebras. It should not be thought of as a generalization of ordinary derived categories over algebras to topological algebras.
\end{remark}
\subsubsection{coHochschild (co-)homology and string topology}
Let $C$ be a (discrete) dg-coalgebra. 
\begin{def1}
	The coHochschild cohomology of $C$ is defined as \[\text{coHH}^n(C) := \Hom_{\Dco(C^{op} \otimes C)}(C, C[n])\]
	and the coHochschild homology of $C$ is 
	 \[\text{coHH}_{n}(C) := H^{-n}(C \square^{R}_{C^e} C)\]
	The cobar construction $\text{Cobar}(C)$ of $C$ is the continuous dual of the bar construction of $C^*$. Alternatively, it can be directly described by a (sum) totalization of a sequence of the form 
	\[ C^{op} \otimes C \to C^{op} \otimes C \otimes C \to C^{op} \otimes C^{\otimes 2} \otimes C \to ...\]
	The coHochschild cocomplex of $C$ is then defined as
	\[ \text{coHC}(C,C) := \Hom_{C^{op} \otimes C}(C,\text{Cobar}(C))[-1].\]
\end{def1}
There is a natural isomorphism $(C^*)^{op} \hat \otimes C^* = (C^{op} \otimes C)^*$. Clearly, the dual of $C$ as a $C^{op} \otimes C$-comodule is again $C^*$ as a pseudocompact $(C^*)^{op} \hat \otimes C^*$-module. Thus, $\text{coHH}^*(C) = \HH^*(C^*)$ and the coHochschild complex of $C$ is naturally isomorphic to the Hochschild complex of $C^*$. In particular, $\text{coHC}(C,C)$ is naturally a dg-Lie algebra.

\begin{proposition}
	\label{PropRivera}
	Let $X$ be a reduced simplicial set whose homotopy category is a groupoid. Then
	\[ \HHa{cont}{*}(C^*(X)) = \mathrm{coHH}_*(C_*(X)) \simeq C_*(L |X|).\]
	\begin{proof}
		The first equality follows from the previous discussion and the second identification is \cite[Theorem 20]{Rivera}.
	\end{proof}
\end{proposition}
\begin{remark}[Hochschild complexes of simplicial sets and string topology]
	\label{RemarkSimplicial}
In \cite{Rivera}, Rivera investigates a certain algebraic model for the free loop space. He starts with a simplicial set $X$ and considers a homotopy equivalent simplicial set $\tilde X$ with the property that every path in $\tilde X$ has a homotopy inverse (all constructions will be independent of the choice of $\tilde X$). The normalized chains $C(\tilde X)$ of $\tilde X$ carry the structure of a so-called \emph{categorical coalgebra}. It is in some sense a multi-object version of a coalgebra; a slightly more general structure is studied in \cite{HolLaza} in which categorical Koszul duality is established. One can then define the coderived category over such categorical coalgebras and define a coHochschild complex. We will show in future work that this notion may also be interpreted in terms of pseudo-compact categories (a definition will be given there). In case $\tilde X$ only has a single $0$-simplex, the coHochschild complex $\text{coHC}_*(C(\tilde X),C(\tilde X))$ of $C(\tilde X)$ as defined above coincides with the definition given in \cite{Rivera} (up to dualizing). As its \emph{homology equals the singular homology of the free loop space of $X$ (Proposition \ref{PropRivera})}, we note that the continuous Hochschild homology of the de-Rham algebra (\ref{TheoremDeRham}) of $X$ (if $X$ is a smooth manifold) is (in general) not equal to the Hochschild homology of the pseudo-compact algebra $C^*(X)$ despite the fact that there is an $A_\infty$-quasi-isomorphism $\Omega(X) \to C^*(X)$ which furthermore induces a so-called \emph{Morita equivalence of the second kind} (as studied in \cite{HC2}), i.e, an equivalence of the dg-categories of finitely generated graded-projective modules. The latter result is known as the \emph{higher Riemann-Hilbert correspondence} first proved in \cite{Block} and then also in \cite{RHHolstein} and generalized in \cite{antweiler}. It thus follows, that Hochschild (co-)homology is \emph{not} an invariant of the (either continuous or discrete) quasi-isomorphism type of the dg-category of finitely generated graded-projective modules. 

 In the case that the simplicial set $X$ is endowed with an appropriately defined notion of an intersection pairing, the coHochschild complex obtains a richer structure providing an algebraic model for the Chas-Sullivan product and Goresky-Hingston coproduct in string topology \cite{RiveraString}. The primary example of such an intersection pairing is induced by Poincaré-duality in case $X$ is a smooth compact oriented manifold. Since the coHochschild complex is isomorphic to the continuous Hochschild complex of the dual pseudo-compact algebra, any extra structure on the coHochschild complex induces additional structure on the continuous Hochschild complex of the pseudo-compact dual.
\end{remark}

\subsubsection{Discrete algebras and Koszul duality}
Let $A$ be a discrete augmented dg-algebra. Its augmented bar construction (not to confuse with the two sided bar construction $\text{Bar}(A)$) $B(A)$ is a coaugmented dg-coalgebra. Then there is an equivalence of triangulated categories \[\text{D}(A) \to \Dco(B(A)),\] see \cite{posi}[Theorem 6.3]. Combining this with the above discussion we see:
\begin{proposition}
	Let $A$ be a discrete augmented dg-algebra and consider its bar construction $B(A)$. Let $C := B(A)^*$ be the dual pseudo-compact dg-algebra (considered as a complete locally convex algebra). Then there is an equivalence of triangulated categories:
	\[ D(A) \to \Dctrpc(C)^{op}\]
	where the latter denotes the full subcategory of $\Dctr(C)$ consisting of pseudo-compact modules.
\end{proposition}
\subsection{Continuous contraderived categories for curved algebras}
So far have considered clct dg-algebras $A$. One of the main advantages of (discrete) contraderived categories is that they may be defined also in the case that the algebra is curved, i.e., if $A$ is a cdg-algebra. 
\begin{def1}
	A \emph{clct cdg-algebra} $A$ is a triple $(A^\#, d, B)$ consisting of a graded clct algebra $A^\#$, a continuous map $d : A \to A$ of degree $1$ and an element $B \in A^2$ such that
	\begin{enumerate}
		\item $d$ is a graded derivation for the product on $A$.
		\item $B$ fulfills the Bianchi identity, $dB = 0$.
		\item $d$ squares to $B$ in the sense that $d^2(x) = Bx - xB$ for all $x \in A$.
	\end{enumerate}
	A (right) clct cdg-module $M$ over $A$ is a graded (right) clct module $M^\#$ over $A^\#$ together endowed with a map $d_M : M \to M$ of degree $1$ fulfilling the graded Leibniz rule and such that $d_M^2(x) = -xB$ for all $x \in M$. As usual, we define the dg-category of right clct cdg-modules $A\text{-Mod}$ (the reader may check that the  differential on Hom-complexes still squares to $0$). Moreover, we may define the contraderived category over $A$ as before.
\end{def1}
\begin{theorem}
	\label{TheoremEnhancement2}
	Let $A$ be a clct cdg-algebra. We denote by $A\text{-Mod}_{gproj}$ the full dg-subcategory of $A\text{-Mod}$ consisting of graded-projective objects. Then the functor 
	\[ H^0(A\text{-Mod}_{gproj}) \to \Dctr(A)\]
	is a triangulated equivalence.
	\begin{proof}
		Note that we have already proven this in the case that $A$ has no curvature $B = 0$, see Theorem \ref{TheoremEnhancement}. In the case that $A$ has curvature, the only difference is that we cannot use the action map $\mu : M \hat \otimes A \to M$ as the first step of a graded-free resolution for $M$ anymore. However, we can use the following construction as used in the proof of \cite[Theorem 3.6]{posi}. For any graded $A^\#$-module $L$, define the cdg-module $G^+(L)$ as follows. Denote the underlying graded vector space of $L$ by $L_{vs}$. Its underyling clct graded vector space is 
		\[G^+(L)_{vs} = L_{vs} \oplus L_{vs}[-1].\] Its elements will be written as $x + dy$ for $x,y \in L$. The $A^\#$ action  on $G^+(L)$ is "twisted":
		\[ \mu(x + dy, a) := xa + (-1)^{|y|}yd(a) + d(ya) = (xa+(-1)^{|y|} yd(a), d(ya)).\]
		$G^+(L)^\#$ is part of a horizontally exact sequence:
		\[ 0 \to L \to G^+(L)^\# \to L[-1] \to 0\]
		of graded $A^\#$-modules.
	 	The differential is given by 
	 	\[ d((x,y)) = d(x+dy) = yB + dx = (yB,x).\]
	 	There is a bijective correspondence between closed degree $0$ maps of cdg-modules $G^+(L) \to M$ and maps of graded $A$-modules $L \to M^{\#}$. This correspondence preserves the class of admissible epimorphisms.
	 	Now choose $L := M \hat \otimes A^\#$. Then the $A$ action on $M$ corresponds to a map $G^+(L) \to M$ which is an admissible epimorphism. Moreover, note that $G^+(L)^\# = L \oplus L[-1]$ in this case now since $L = M \hat \otimes A^\#$ is graded-free. So $G^+(L) \to M$ yields the first step of a graded-free resolution of $M$. We may thus proceed iteratively and totalize at the end.
	\end{proof}
\end{theorem}
\subsubsection{Generalized resolutions and the curved bar construction}
In the following, we would like to have an analogue of the Bar construction. For curved algebras however, we get additional differentials to the usual horizontal and vertical differentials (an extra term decreasing the horizontal degree by $1$ and increasing the vertical by $2$). For such kind of "double complexes", there is a generalization of Lemma \ref{LemmaTot}. Lemma \ref{LemmaTot} is recovered when setting $d_1 = d_2 = ... = 0 $. Moreover, we will only use the result of the Lemma in case $d_2 = d_3 = ... = 0$.
\begin{lemma}
	\label{lemmaGenResolution}
	Let $A$ be a clct cdg-algebra. Let $M^i$, $i \leq 0$ be a collection of graded $A^\#$-modules together with the following data:
	\begin{enumerate}
		\item Continuous $\mathbb C$-linear maps $d_v : M^i \to M^i$ of degree $1$ which fulfill the Leibniz rule with respect to the differential $d$ on $A$.
		\item Continuous $A^\#$-linear maps $d_k : M^i \to M^{i-k}$ for $k \in \mathbb N_{\geq 1} \cup \{-1\}$ of degree $k+1$ such that $d = d_{-1} + d_v + d_1 + d_2 + ...$ defines a differential on the product totalization $\text{Tot}^{\Pi}(M)$ of $M^*$ (i.e, $d^2(x) = -xB$, note that $d$ fulfills the Leibniz rule by construction).
	\end{enumerate}
	In that case $\text{Tot}^{\Pi}(M)$ is naturally a cdg $A$-module. Note that $(d_{-1})^2 = 0$ for degree reasons. If $d_{-1}$ admits a continuous $\mathbb C$-linear contracting homotopy, then $\text{Tot}^{\Pi}(M)$ is contraacyclic.
	\begin{proof}
		In view of Theorem \ref{TheoremEnhancement2}, it suffices to prove that for each cdg-module $N$ which is graded-free, we have that 
		\[ \Hom^*_{A}(N, \text{Tot}^{\Pi}(M)) \simeq 0.\]
		Thus, write $N^\# = V \hat \otimes A^\#$. Observe that all the differentials $d_i$ and $d_v$ do not decrease vertical degrees. Write the underlying graded space of $\Hom^*_{A}(N, \text{Tot}^{\Pi}(M))$ as
		\[ \Hom^*_{A}(N, \text{Tot}^{\Pi}(M))^\# = \Hom^*_{\mathbb C}(V, \text{Tot}^{\Pi}(M))^\#\]
		The latter space is naturally filtered by the difference of the vertical degrees of the $M^i$ and $V$. It is easy to see that this filtration is compatible with the induced differential (as all differentials only increase the vertical degree). The filtration is degreewise bounded above and complete. The associated spectral sequence thus converges to the cohomology of $ \Hom^*_{A}(N, \text{Tot}^{\Pi}(M))$. The first page of the spectral sequence is given by the $d_{-1}$ cohomology\[H^*(\Hom^*_{\mathbb C}(V, (\text{Tot}^{\Pi}(M),d_{-1}))).\] Thus, if $d_{-1}$ admits a continuous $\mathbb C$-linear contracting homotopy, then the first page is trivial and we obtain 
		\[ \Hom^*_{A}(N, \text{Tot}^{\Pi}(M)) \simeq 0.\]
	\end{proof}
\end{lemma}
\begin{def1}
	\label{GenRes}
	A resolution as in the statement of the previous Lemma will be called a \emph{generalized resolution}.
\end{def1}
\begin{def1}
	\label{DefResolutionCurved}
	Here we will give a detailed description of the two-sided bar construction. Let $(A, d, B)$ be a clct cdg-algebra. It will be a resolution of $A$ as a $A^{op} \hat \otimes A$-module in the generalized sense, Definition \ref{GenRes}. We first define graded $(A^\#)^{op} \hat \otimes A^{\#}$-modules $\text{Bar}^{-n}(A)$ for $n \geq 0$:
	\[ \text{Bar}^{-n}(A) := A^{op} \hat \otimes A^{\hat \otimes n} \hat \otimes A\]
	each of which is endowed with a map $d_v$ given on pure tensors by
	\[ d_v(a_0 \hat \otimes ...  \hat \otimes a_{n+1}) := (-1)^{n} \sum_{i = 0,...,n+1} (-1)^{|a_0| + ... + |a_{i-1}|} a_0 \hat \otimes ... \hat \otimes d(a_i) \hat \otimes ... \hat \otimes a_{n+1}.\]
	It does not square to $0$ (unless $B = 0$) but it fulfills the Leibniz rule. Furthermore, we define
	\begin{nalign} &d_{-1} : \text{Bar}^{-n} \to \text{Bar}^{-n+1} \\ &d_{-1}(a_0 \hat \otimes ...  \hat \otimes a_{n+1}) := \sum_{i = 0,...,n} (-1)^i a_0 \hat \otimes ... \hat \otimes a_i a_{i+1} \hat \otimes ...  \hat\otimes a_{n+1}
	\end{nalign}
	and the extra piece $d_1$ involving the curvature $B$:
	\begin{nalign} &d_{1} : \text{Bar}^{-n} \to \text{Bar}^{-n-1} \\ &d_{1}(a_0 \hat \otimes ...  \hat \otimes a_{n+1}) := \sum_{i = 1,...,n+1} (-1)^i a_0 \hat \otimes ... \hat \otimes B \hat \otimes a_i \hat \otimes ...  \hat\otimes a_{n+1}
	\end{nalign}
	these maps fulfill the following relations:
	\[ (d_1)^2 = 0, (d_{-1})^2 = 0, d_v d_1 + d_1 d_v = 0, ((d_v)^2 + d_{-1} d_1 + d_1 d_{-1})(x) = Bx - xB.\]
	We set $\text{Bar}(A) := (\prod_{n \geq 0} \text{Bar}^{-n}(A)[-n], d = d_{-1} + d_v + d_1)$
	which defines an $A^{op} \hat \otimes A$-module. 
	Furthermore, multiplication $A^{op} \hat \otimes A\to A$ induces an augmentation $\text{Bar}^*(A) \to A$. The map $d_{-1}$ (together with the augmentation) admits a continuous $\mathbb C$-linear contracting homotopy:
	\[h(a_0 \hat \otimes ... \hat \otimes a_{n+1}) := 1 \hat \otimes a_0 \hat \otimes ... \hat \otimes a_{n+1}.\]
	Thus, in view of Lemma \ref{lemmaGenResolution}, we have that $\text{Bar}(A)$ is a graded-free replacement for $A$ as a cdg $A^{op} \hat \otimes A$-module. In a similar spirit, we may define the reduced bar construction $\text{Bar}_{red}(A)$ and obtain the same statements.
\end{def1}
\begin{def1}
	Let $A$ be a clct cdg-algebra. We define the \emph{ continuous Hochschild complex of $A$ with values in $M$} as
	\[ \HC(A, M) := \Hom_{A^e}(\text{Bar}_{red}(A), M)[-1]\]
	Its cohomology will be called the \emph{Hochschild cohomology of $A$ with values in $M$} and by the previous results it is isomorphic to
	\[ \Hom_{\Dctr(A^e)}(A[n],M).\] 
	The results of Theorem \ref{IsDgLie} carry over to this case so $\HC(A, A)$ is a dg-Lie complex in a canonical way.
\end{def1}

\subsection{Comparison with the discrete case}
Our construction begs to be compared with Positselski's original construction of contraderived categories over discrete algebras. There is a fully-faithful functor from discrete vector spaces to complete locally convex vector spaces (it is left-adjoint to the forgetful functor), see \cite[p. 111]{Jarchow} and Lemma \ref{DiscreteLocallyConvex}.
\begin{lemma}
	Let $A$ be a discrete graded vector space. Viewing $A$ as a graded clct space, if $A$ carries the structure of a clct cdg-algebra (as in Definition \ref{DefAlgebras}) then it canonically carries the structure of a discrete cdg-algebra in the usual sense. The converse is false. However, if $A$ has a countable basis, then the structure of a discrete cdg-algebra on $A$ gives rise to a clct cdg-algebra if and only if for each $k$ there are only finitely many $i \in \mathbb Z$ such that the multiplication $A^{k-i} \otimes A^{i} \to A^k$ is non-zero.
	\begin{proof}
		The point is the following identity for discrete spaces
		\[ \Hom_{\text{continuous}}(\prod V_i, W) = \oplus \Hom_{\text{continuous}}(V_i, W) = \oplus \Hom_{\text{discrete}}(V_i, W) \]
		because $W$ carries a Hausdorff seminorm by existence of bases (see Lemma \ref{LemmaProdPlus}). Furthermore, 
		\[ \oplus \Hom_{\text{discrete}}(V_i, W)  \subseteq \prod \Hom_{\text{discrete}}(V_i, W).\]
		For the second claim it is required that $A$ has a countable basis so that the discrete tensor product agrees with the projective one, see Lemma \ref{InductiveLimits}.
	\end{proof}
\end{lemma}
Now let $A$ be a discrete space endowed with the structure of a clct cdg-algebra. By the previous lemma, we may also view it as a discrete cdg-algebra so the contraderived category of $A$ in the (discrete) sense of Positselski is defined, let us denote it by $\D{ctr,P}{}(A)$. On the other hand, we may use our definition of the contraderived category $\Dctr(A)$. The forgetful functor from complete locally convex spaces to discrete spaces is a right-adjoint and hence preserves products. It follows that there is a natural comparison functor $\Dctr(A) \to \text{D}^{\text{ctr,P}}(A)$. Of course it is not expected that this functor is an equivalence, but one might for example try to restrict to the subcategory of modules over $A$ that are discrete, let us denote this full subcategory by $\D{ctr}{disc}(A)  \subseteq \Dctr(A)$.

We remind the reader of the construction of \emph{Ext of the second kind}, following \cite{Polishchuk}. Take a discrete cdg-algebra $A$ and two cdg-modules $M$ and $N$. Construct a horizontally exact resolution $... \to F_2 \to F_1 \to F_0 \to M$ of cdg-modules by graded-projectives $F_i$ (in the discrete sense). Then define
\[ \Exta{II, PP}{}(M,N) := \text{Tot}^{\oplus} \Hom^*_{A}(F_*, M). \] 
It is a rather unfortunate fact that $\Exta{II, PP}{}(M,N)$ does \emph{not} compute the contraderived Hom of $M$, $N$ (it is not even invariant under contra-isomorphisms in $M$). Our theory fixes this issue (one has to use a product totalization inside of Hom in general, but for a discrete module $M$, the totalization can be pulled out). It is now tempting to conclude that the difference between $\D{ctr}{disc}(A)$ and  $\D{ctr,P}{}(A)$ lies in the fact that in the former one computes morphisms as Ext of the second kind as defined by Polishchuk-Positselski. However, this is not quite the case since a discrete graded-projective module might not be graded-projective in the continuous sense.
\begin{proposition}
	\label{LemmaComparison}
	Suppose $A$ has a countable basis. Let $M$ and $N$ be clct cdg-modules over $A$ which carry the discrete topology and such that $M$ has a countable basis. Then
	\[ \Hom_{\Dctr(A)}(M,N) \cong \Exta{II, PP}{}(M,N).\]
	\begin{proof}
		Both sides may be computed in terms of bar resolutions. On the left hand side we use the projective tensor product while on the right side the discrete one is employed. As $M$ and $A$ have countable bases, these tensor products agree by \cite[Corollary 15.5.4]{Jarchow}.
	\end{proof}
\end{proposition}
Note also that if $M$ is finitely-generated graded-projective, then $\Hom$ in $\Dctr(A)$ agrees with that in $\text{D}^{\text{ctr,P}}(A)$ and also with Ext of the second kind.
\begin{corollary}
	\label{CorComparison}
	Let $A$ be a clct cdg-algebra endowed with the discrete topology and suppose that $A$ has a countable basis. Then 
	\[ \HHa{PP}{}(A,A) \cong \Hom_{\Dctr(A^e)}(A,A) =: \HHcont(A,A).\]
\end{corollary}

\section{Continuous formality theorems for manifolds}

\subsection{Smooth functions}
It is well known that the continuous Hochschild cohomology of the Fréchet algebra of smooth functions is given by $\Lambda \Gamma(TM)$, the space of polyvectorfields \cite{Pflaum}. Moreover, the usual definition of continuous Hochschild cohomology agrees with ours as it is concentrated in degree $0$, see Definition \ref{DefHochschild}. It is easy to show the above result in the case that $M$ is an open subset of $\mathbb R^n$. However, the globalisation techniques in the literature are rather involved, use tubular neighbourhood theorems and exponential maps associated to Riemannian metrics. In our framework, globalisation is entirely straightforward.
\begin{def1}
	Consider the manifold $\mathbb R^n$ and its Fréchet algebra of smooth $\mathbb K$-valued functions $C^\infty(\mathbb R^n)$. We consider the $C^\infty(\mathbb R^n \times \mathbb R^n)$-modules
	\[ F^{std,-k} := \Lambda^k {\mathbb K^n} \otimes C^\infty(\mathbb R^n \times \mathbb R^n)\]
	for every $k \geq 0$. The standard coordinates on $\mathbb R^n \times \mathbb R^n$ are denoted by \\ $x_1,...,x_n,x'_1,...,x'_n$. Let us also define a new set of coordinates $t_i := x_{i} - x'_{i}$ and $s_i := \frac{x_i+x'_i}{2}$. Moreover, the basis of $\mathbb K^n$ is denoted by $dt_1^s,...,dt_n^s$. We define a horizontal differential 
	\[ d_h(\alpha \otimes f) = \sum_{i = 1}^{2n} \iota_{\partial_{t^s_i}} \alpha \otimes t_i f .\]
	Moreover, there is an obvious map $F^{std,0} \to C^\infty(\mathbb R^n)$ of $C^\infty(\mathbb R^n \times \mathbb R^n)$-modules. Clearly, each $F^{std,-k}$ is free and moreover we have:
\end{def1}
\begin{lemma}
	\label{RealRes}
	The complex \[ ... \to F^{std,-1} \to F^{std,0} \to C^\infty(\mathbb R^n)\] 
	admits a continuous $\mathbb K$-linear contracting homotopy.
	\begin{proof}
		We construct the homotopy in the case that $n = 1$. The general case follows by observing that the sequence for higher $n$ arises as the $n$-fold projective tensor product of the $1$-dimensional case. We expand any function $f$ on  $\mathbb R^{n} \times \mathbb R^{n}$ as 
		\[ f(t,s) = f(0,s) + t \tilde f(t,s)\]
		where we note that $f \mapsto \tilde f$ is continuous.
		The homotopy for
		\[
		\begin{tikzcd}
			0 \arrow[r] & C^{\infty}(\mathbb R^2) \arrow[r, "\cdot t"'] & C^{\infty}(\mathbb R^2) \arrow[r, "\text{ev}_{t = 0}"'] \arrow[l, "h_1"', dashed, bend right] & C^{\infty}(\mathbb R)  \arrow[r] \arrow[l, "h_2"', dashed, bend right] & 0
		\end{tikzcd}
		\]
		 is then given by 
		 \begin{align*}
		 	h_1(f) = \tilde f \text{ and } h_2(g)(t,s) = g(0,s).
		 \end{align*}
	\end{proof}
\end{lemma}
\begin{theorem}
	\label{HochschildFunctions}
	Let $M$ be a smooth manifold and denote by $C^{\infty}(M)$ the Fréchet algebra of smooth functions on $M$. Then the HKR-map 
	\begin{nalign}\label{CohSmooth} \Lambda^* \Gamma(TM)[-1] \to \HC(C^{\infty}(M),C^{\infty}(M)) \end{nalign}
	is a quasi-isomorphism. Moreover, the above map can be extended to an $L_\infty$-quasi isomorphism
	\begin{nalign}\Lambda^* \Gamma(TM)[-1] \to \HC(C^{\infty}(M),C^{\infty}(M)) \end{nalign}
	of dg-Lie algebras.
	\begin{proof}
		We begin by verifying the first claim.
		Standard Koszul resolutions establish the result locally as we have seen in the previous Lemma. The assignment
		\[ U \mapsto \HC(C^{\infty}(M), C^{\infty}(U))\]
		defines a complex of sheaves of $C^{\infty}(M)$-modules on $M$ (since $U \mapsto C^{\infty}(U)$ is a Fréchet sheaf). As a consequence of the following Lemma, we have the identity
		\[ C^{\infty}(U) = C^{\infty}(M) \hat \otimes_{C^\infty(M \times M)} C^{\infty}(U \times U)\]
		and by Example \ref{ExRes}, we have the functor
		\[ f_* = (-) \hat \otimes_{C^{\infty}(M \times M)} C^{\infty}(U \times U) : \Dctr(C^{\infty}(M \times M)) \to \Dctr(C^{\infty}(U \times U))\]
		which sends $C^\infty(M)$ to $C^{\infty}(U)$. It follows that the natural map
		\begin{nalign}\label{LocalReduction} \text{HC}_{C^\infty(M \times M)}^{cont}(C^{\infty}(M), C^{\infty}(U)) \simeq \text{HC}_{C^\infty(U \times U)}^{cont}(C^{\infty}(U), C^{\infty}(U)) \end{nalign}
		is a quasi-isomorphism for every open $U \subseteq M$ by Lemma \ref{Adjunction}. We can extend the map \ref{CohSmooth} to a map of sheaves of $C^{\infty}(M)$-modules which is a quasi-isomorphism of (bounded below) complexes of sheaves by the local result and equation \ref{LocalReduction} (there is an obvious commutative diagram). As we are working over the fine sheaf $C^{\infty}(M)$, the quasi-isomorphism induces a quasi-isomorphism on global sections.
		
		Let us now turn to the second statement regarding the extension to an $L_\infty$-quasi isomorphism. Kontsevich \cite{Kontsevich} showed the analogous statement for the Hochschild complex of polydifferential cochains which we denote by $\text{HC}_{\text{PD}}(C^\infty(M))$. There is a natural inclusion of dg-Lie algebras  $\text{HC}_{\text{PD}}(C^\infty(M)) \to \HC(C^\infty(M))$. Composing this dg-Lie algebra map with Kontsevich's $L_\infty$-quasi-isomorphism 
		\begin{nalign} \Lambda^* \Gamma(TM)[-1] \to \text{HC}_{\text{PD}}(C^{\infty}(M),C^{\infty}(M)) \end{nalign}
		yields the desired $L_\infty$-quasi isomorphism (using that an $L_\infty$-map is a quasi-isomorphims if and only its zeroth component is a quasi-isomorphism of complexes).
	\end{proof}
\end{theorem}
\begin{lemma}
	\label{denseImage}
	Let $p : A \to B$ be a map of clct dg-algebras which has dense range.
	Then the natural map 
	\[ s : A \hat \otimes_{A^e} B^e \to B\]
	is an isomorphism. 
	\begin{proof}
		There is another continuous map $i : B \to A \hat \otimes_{A^e} B^e, i(b) = 1 \otimes (b \otimes 1)$. Clearly, $s \circ i = \text{id}_B$. The result follows if we can show that $i$ has dense range. This can be deduced by considering the following commutative diagram
		\[\begin{tikzcd}
			A \hat \otimes_{A^e} A^e \arrow[r, "\text{id} \otimes p \otimes p"]         & A \hat \otimes_{A^e} B^e \\
			A \arrow[r, "p"] \arrow[u, no head, equal] & B \arrow[u, "i"]            
		\end{tikzcd}\]
		and using the density of the upper row.
	\end{proof}
\end{lemma}

\subsection{Dolbeault algebra}
Let $X$ be a complex manifold and consider its Dolbeault Algebra $\mathcal A(X) := (\Omega^{0,*}(X), \bar \partial)$, which we view as a Fréchet dg-algebra. We can reduce the computation of the contraderived Hochschild cohomology of $\mathcal A(X)$ in the same fashion as in Theorem \ref{HochschildFunctions}, to a local one. The local resolutions are a bit more involved, so we will construct them here. 

\begin{lemma}
	Consider the complex manifold $\mathbb C^n$ and the point $0 \in \mathbb C^n$. Consider $\mathbb C$ as a $\mathcal A(\mathbb C^n)$ module with the action $1.f := 1 \cdot f(0)$. There is a graded-free module of the form $R = (\Lambda^{-*} \langle d\bar z_1^f,...,d\bar z_n^f, d z_1^f, ..., d z_n^f \rangle \otimes \mathbb C[[\bar z_1^f,...,\bar z_n^f]] \hat \otimes \mathcal A(\mathbb C^n), d_{\mathcal A} + \mathbb E^0_h + \mathbb E^1_h + \mathbb E_v)$ where
	\begin{enumerate}
		\item $|dz_i^f| = |d \bar z_i^f| = -1$, $|\bar z_i^f| = 0$,
		\item $d_{\mathcal A} = 1 \otimes 1 \otimes \bar \partial $ is the usual Dolbeault differential acting on the right-most factor,
		\item $\mathbb E_h$ is given by
		\[ \mathbb E^0_h(\omega \otimes f \otimes 1) = \sum_{i = 1}^n \iota_{\bar \partial_i} \omega \otimes f \otimes \bar z_i + \iota_{ \partial_i} \omega \otimes f \otimes  z_i  \]
		and 
		\[ \mathbb E^1_h(\omega \otimes f \otimes 1)  = \sum_{i = 1}^n (-1)^{|\omega|} \omega \otimes \bar \partial_i f \otimes d \bar z_i\]
		\item $\mathbb E_v$ is defined as
		\[\mathbb E_v(\omega \otimes f \otimes 1) = - \sum_{i = 1}^n  \iota_{\bar \partial_i} \omega \otimes \bar z_i^f f \otimes 1\]
	\end{enumerate}
	and a contraderived isomorphism $R \to \mathbb C$.
	\begin{proof}
		One easily verifies the relations
		\[ (\mathbb E^0_h)^2 = (\mathbb E^1_h)^2 = (\mathbb E_v)^2 = [d_{\mathcal A},\mathbb E_h^1] = [d_{\mathcal A},\mathbb E_v] = [\mathbb E_h^0,\mathbb E_v] = [\mathbb E_h^0,\mathbb E_h^1] = 0\]
		and 
		\[ [\mathbb E_h^1,\mathbb E_v] + [d_{\mathcal A},\mathbb E_h^0] = 0\]
		and thus we can view $R$ as the product totalization of a bounded above sequence of $\mathcal A(\mathbb C^n)$-modules for which $\mathbb E_h = \mathbb E^0_h + \mathbb E^1_h$ is the horizontal differential. Moreover, there is an augmentation $R \to \mathbb C$ induced by $\text{ev}_0 : \mathcal A(\mathbb C^n) \to \mathbb C$.
		To see that the double complex is horizontally exact, we may forget about $\mathbb E_v$. The resulting resolution is a resolution of $\mathbb C \otimes \mathbb C$ over the graded algebra $C^{\infty}(\mathbb C^n) \hat \otimes \Lambda^* \mathbb C^n$. More precisely, it is the tensor product of the resolution of $\mathbb C$ as a $C^{\infty}(\mathbb C^n)$-module provided by Lemma \ref{RealRes} with the standard polynomial Koszul resolution $\mathbb C[[\bar z_1^f,...,\bar z_n^f]] \hat \otimes \Lambda^* \mathbb C^n$ of $\mathbb C$ as a module over $\Lambda^* \mathbb C^n$ ($\bar z_i^f$ has bidegree $(-1,1)$ in this resolution). An explicit continuous contracting homotopy for the latter resolution is given the formula
		\[ 	H(f \otimes \omega) := \frac{(-1)^{p+1}}{p+q} \sum_{i=1}^n f \cdot z^i \otimes \iota_{\partial_i} \omega\]
		for $f$ a homogeneous polynomial of degree $q$ and $p$ is the degree of $\omega$. 
	\end{proof}
\end{lemma}

\begin{proposition}
	\label{LocalDolbeault}
	Let $U \subseteq \mathbb C^n$ be an open subset. Then the HKR-map
	$\Lambda T^{1,0} U \otimes \mathcal A(U)[-1] \to \HC(\mathcal A(U))$
	is a quasi-isomorphism.
	\begin{proof}
		We start with the case that $U = \mathbb C^n$. By an affine transformation, we can consider $\mathbb C^n \to \mathbb C^n \times \mathbb C^n, z \mapsto (0,z)$ instead of the diagonal embedding. Then $\mathcal A(\mathbb C^n) = \mathbb C \otimes \mathcal A(\mathbb C^n)$ has the resolution $R \hat \otimes \mathcal A(\mathbb C^n)$ as constructed in the previous Lemma. Let $R_{d_h = 0}$ be the module obtained by setting the horizontal differential of $R$ to $0$ (i.e., replace $d_{\mathcal A} + \mathbb E^0_h + \mathbb E^1_h + \mathbb E_v$ by $d_{\mathcal A} + \mathbb E_v$). In the complex
		\[ \text{Hom}^*_{\mathcal A(\mathbb C^{2n})}(R \hat \otimes \mathcal A(\mathbb C^n), \mathbb C \hat \otimes \mathcal A(\mathbb C^n))\]
		the horizontal differential $\mathbb E^0_h + \mathbb E^1_h$ of $R$ acts as zero (since $\bar z_i, z_i$ and $d \bar z_i$ act as zero on $\mathbb C$) and so 
		\[ \text{Hom}^*_{\mathcal A(\mathbb C^{2n})}(R \hat \otimes \mathcal A(\mathbb C^n), \mathbb C \hat \otimes \mathcal A(\mathbb C^n)) = \text{Hom}^*_{\mathcal A(\mathbb C^{2n})}(R_{d_h = 0} \hat \otimes \mathcal A(\mathbb C^n), \mathbb C \hat \otimes \mathcal A(\mathbb C^n)).\]
		Now consider the submodule of $R_{d_h = 0}$ given by $R^{red} := \Lambda^{-*} \langle d z_1^f, ..., d z_n^f \rangle \otimes \mathcal A(\mathbb C^n)$. There are $\mathcal A(\mathbb C^n)$-linear continuous homotopy equivalences $R^{red}_i \simeq R_{d_h = 0, i}$ for each $i$ where $i$ is the horizontal degree. As the horizontal differential is equal to $0$, these equivalences assemble to an equivalence in the contraderived category and so 
		\[ \text{Hom}^*_{\mathcal A(\mathbb C^{2n})}(R_{d_h = 0} \hat \otimes \mathcal A(\mathbb C^n), \mathbb C \hat \otimes \mathcal A(\mathbb C^n)) = \text{Hom}^*_{\mathcal A(\mathbb C^{2n})}(R^{red} \hat \otimes \mathcal A(\mathbb C^n), \mathbb C \hat \otimes \mathcal A(\mathbb C^n))\]
		and the latter complex is equal to 
		\[ \Lambda T^{1,0} \mathbb C^n \otimes \mathcal A(\mathbb C^n).\]
		
		Now let $U \subset \mathbb C^n$ be an arbitrary open subset. The resolution $R \hat \otimes \mathcal A(\mathbb C^n)$ can be transformed back to a resolution of the diagonal module $\mathcal A(\mathbb C^n)$. By the projectivity of $\mathcal A(U \times U)$ over $\mathcal A(\mathbb C^{2n})$, all complexes considered before can be tensored with $\mathcal A(U \times U)$ and then one may argue as in the case of $\mathbb C^n$.
	\end{proof}	
\end{proposition}

\begin{theorem}
	\label{CtrDolbeault}
	Let $X$ be a complex manifold.
	The HKR-map 
	\[ \Lambda T^{1,0} X \otimes \mathcal A(X)[-1] \to \HC(\mathcal A(X))\]
	is a quasi-isomorphism. Moreover, this map can be extended to an $L_\infty$-quasi-isomorphism between $\Lambda T^{1,0} X \otimes \mathcal A(X)[-1]$ and $\HC(\mathcal A(X))$. In particular, 
	\[ \HH_{\text{cont}}(\mathcal A(X)) = H(X, \Lambda T^{1,0}).\]
\begin{proof}
		In the previous Proposition, we have established that \[\Lambda T^{1,0} U \otimes \mathcal A(U) \to \HC(\mathcal A(U))\] is a quasi-isomorphism if $U$ is biholomorphic to an open subset of $\mathbb C^n$. As in the proof of Theorem \ref{HochschildFunctions}, we obtain that it must a quasi-isomorphism globally. Actually, we need to be a little bit more careful as the Hochschild complex is now also unbounded in the negative direction but the statement remains valid because of Lemma \ref{LemmaGlobalisation}.
		
		The second statement can be proven as in Theorem \ref{HochschildFunctions} using the analogous holomorphic formality results, see \cite[Corollary 5.3]{Damien}.
	\end{proof}
\end{theorem}
\begin{remark}
	This theorem is yet another manifestation of the crucial rôle of the Dolbeault algebra in the theory of complex manifolds. It has been proven by Block \cite{Block}, that the category of cohesive modules over $\mathcal A(X)$ is equivalent to the bounded derived category complexes of sheaves with coherent cohomology in case $X$ is compact, see also \cite{Bondal}. This has been generalized to non-compact manifolds in \cite{RHHolstein} and furthermore in the presence of a (topologically trivial) gerbe in \cite{antweiler}. We remark that the dg-category of cohesive modules is equivalent to the dg-category of finitely generated $\mathcal A(X)$-modules that are graded-projective. As we have seen, this is a full subcategory of $\Dctr(\mathcal A(X))$. It would be interesting to study $\Dctr(\mathcal A(X))$ as a candidate for a  category of analytic quasi-coherent sheaves on $X$, see subsection \ref{subsectionSheaves} for some context.
\end{remark}
\begin{lemma}
	\label{LemmaGlobalisation}
	Let $X$ be a Hausdorff topological space of finite covering dimension and $F$ a (possibly unbounded) complex of soft sheaves on it. Suppose for each $x$ there exists a neighbourhood $U$ such that for every open subset $V \subseteq U$ we have that $F^*(V)$ is an acyclic complex (of abelian groups). Then $F^*(X)$ is acyclic.
	\begin{proof}[Proof Sketch]
		By the assumptions on $X$ one can show that there exists a cover $(U_i)_{i \in I}$ of $X$ with $I$ finite such that $F^*(U)$ is acyclic for each finite intersection $U$ of the $U_i$. We may then regard the alternating \v Cech complex of $F$ with respect to the cover as a horizontally bounded double complex. Thus, there are 2 associated spectral sequences converging to the cohomology of its totalization. The first spectral sequence collapses to $0$ on the first page. The second one collapes on the second page to $H^*(F^*(X))$ (using softness).
	\end{proof}
\end{lemma}
\subsubsection{Twisted Dolbeault algebra and generalized complex manifolds}
In this subsection we will generalise Theorem \ref{CtrDolbeault} in the presence of a topologically trivial gerbe, i.e., a $\bar \partial$-closed $(0,2)$-form $B$ on $X$. The resulting geometric space can naturally be interpreted in the framework of \emph{generalised complex geometry} introduced by Hitchin and Gualtieri, for background we refer to \cite{gualtieri}.
\begin{theorem}
	\label{CtrDolbeault2}
	Let $X$ be a complex manifold and $B \in \mathcal A^{0,2}(X)$ a $\bar \partial$-closed $(0,2)$-form on $X$. And consider the cdg-algebra $\mathcal A_B := (\mathcal A(X), \bar \partial, B)$. To $B$ we can associate a generalized complex manifold $X_B$. There is a quasi-isomorphism of complexes:
	\[ \Gamma(\Lambda^* L_B)[-1] \to \HC(\mathcal A_B)\]
	where the left side denotes the deformation dg-Lie complex associated to $X_B$.
	In case $X$ is a compact Kähler manifold, this map can be extended to an $L_\infty$-quasi-isomorphism of dg-Lie complexes. Moreover, in that case, both dg-Lie complexes are also $L_\infty$-quasi-isomorphic to their untwisted versions.~
	\begin{proof}
		The proof is mainly an analysis on what effects $B$ has on both sides. Following \cite{gualtieri}, the Lie-algebroid $L_B$ is given by 
		\[ L_B = \{ X + \iota_X(\bar B) | X \in T^{1,0}\} \oplus (T^{0,1})^{*}\]
		where $\bar B \in \mathcal A^{2,0}(X)$ is the complex conjugate to $B$.
		It is endowed with the Dorfman bracket which is continued to $\Gamma(\Lambda^* L_B)$ to yield a graded Lie-bracket on that space. The differential on $\Gamma(\Lambda^* L_B)$ is induced by the Chevalley-Eilenberg differential for the Lie-algebroid $\bar L_B$ (using $\bar L_B = L_B^*$ by means of the scalar product). There is an isomorphism of vector bundles 
		\[ \phi : L = T^{1,0} \oplus (T^{0,1})^* \to L_B, \phi(X + \xi) := X + \iota_X(\bar B) + \xi\]
		(and similar for their conjugates).
		Note that $\phi$ is in general not a morphism of Lie-algebroids. Using $\phi$, we may transport the Lie-algebroid structure to $L = T^{1,0} \oplus (T^{0,1})^*$. Denote the resulting Lie-bracket on $L$ by $[-,-]_B$ while the untwisted one is denoted by $[-,-]$. They are related by the following formula:
		\begin{nalign} \label{EqTwist} [X + \xi, Y + \eta]_B = [X + \xi, Y + \eta] + \bar \partial \bar B(X,Y,-).
			\end{nalign}
		The induced B-twisted differential on $\Gamma(\Lambda^* L)[-1]$ has the form
		\[ d_B(a) = d(a) + [B,a]\]
		where $d$ denotes the untwisted differential. On the right hand side, the by $B$-twisted differential on $\HC(\mathcal A_B)$ takes a similar form:
		\[ d^{HH}_B(\mu) = d^{HH}(\mu) + [B,\mu]\]
		which can be seen from Definition \ref{DefResolutionCurved}.
		The key fact which we are going to use is that the usual $HKR$-map (of untwisted dg-Lie algebras)
		\[ l_0 : \Gamma(\Lambda^* L)[-1] \to \HC(\mathcal A)\]
		fulfills $l_0([a,b]) = [l_0(a),l_0(b)]$ whenever either $a$ or $b$ is contained in \[\mathcal A^{0,*} \subseteq \Gamma(\Lambda^* L)\] which can be seen from a local computation. Thus, it is also true that $l_0$ defines a map of $B$-twisted complexes 
		\[ (\Gamma(\Lambda^* L), d_B)[-1] \to \HC(\mathcal A_B).\]
		We remark that it may be possible that the same holds true for the entire $L_\infty$-morphism, but we will not attempt to do this here. Now, to check that it is a quasi-isomorphism we can reduce it to a local computation by the same technique as used in the proof of Theorem \ref{CtrDolbeault}. Locally, $B$ defines a trivial cohomology class, and one can use one of its primitives to untwist both complexes (if $B = \bar \partial A$, then $\mu \mapsto \mu + [A,\mu]$ defines gauge transformations on both sides). Thus, 
		\begin{nalign}
			\label{EqCompare}
		(\Gamma(\Lambda^* L), d_B)[-1] \to \HC(\mathcal A_B)
		\end{nalign}
		is a quasi-isomorphism globally.
		
		We now let $X$ be a compact Kähler. In that case, in virtue of the Hodge decomposition, we can find a $(0,2)$-form $\tilde B$ which defines the same Dolbeault cohomology class as $B$ but is $d$-closed (i.e, it is $\partial$-closed and $\bar \partial$-closed). As seen above, passing to a different representative yields equivalent dg-Lie complexes on both sides of \ref{EqCompare}. Observe that the twist on the left side of \ref{EqCompare} by $\tilde B$ is $0$ , as the twist only depends on $\bar \partial \bar{\tilde B} = 0$ (see equation \ref{EqTwist}). On the other side, twisting by $\tilde B$ does change the differential but we can argue in the following manner. There is another dg-Lie complex consisting of $\mathcal O_X$ polydifferential operators and globalized using $\mathcal A(X)$:
		\[ \Gamma(\text{PolyDiff}(\mathcal O_X) \otimes_{\mathcal O_X} \mathcal A)\]
		which is the dg-Lie complex considered in \cite{Damien}. There is a natural inclusion of dg-Lie complexes
		\[ \Gamma(\text{PolyDiff}(\mathcal O_X) \otimes_{\mathcal O_X} \mathcal A) \to \HC(\mathcal A)\]
		which, as seen before, is a quasi-isomorphism. But the same map also defines a morphism
		\[ \Gamma(\text{PolyDiff}(\mathcal O_X) \otimes_{\mathcal O_X} \mathcal A) \to \HC(\mathcal A_{\tilde B})\] 
		because every cochain $\mu$ in the image of this map has the property $[\tilde B, \mu] = 0$ (using that $d \tilde B = 0$). Moreover, we can argue that it is a quasi-isomorphism by sheaf-theoretic arguments as before.
		Thus, there exists an $L_\infty$-quasi-isomorphism between $\HC(\mathcal A)$ and $\HC(\mathcal A_B)$.

	\end{proof}
\end{theorem}
\begin{remark}
	We expect that 
	\[ \Gamma(\Lambda^* L_B) \to \HC(\mathcal A_B)\] 
	can be extended to an $L_\infty$-quasi-isomorphism for any complex manifold $X$. A proof of this statement seems to require a careful inspection of the techniques used in \cite{Damien} and \cite{Kontsevich} which is beyond the scope of this paper.
\end{remark}
\subsection{de Rham algebra}
Let $M$ be a smooth manifold and denote by $\Omega(M)$ its de Rham algebra (real or complex valued). For any clct dg-algebra $A$ there is an inclusion $i : Z(A[-1]) \to \text{HC}_{cont}(A)$ of dg-Lie algebras where $Z(-)$ denotes the graded center endowed with the trivial Lie-bracket. In the case of the de-Rham algebra this map is a quasi-isomorphism:
\begin{theorem}
	\label{TheoremDeRham}
	Let $M$ be a smooth manifold, then the inclusion
	\[i : \Omega(M)[-1] \to \HC(\Omega(M))\]
	is a quasi-isomorphism of dg-Lie algebras where the left side is endowed with the trivial bracket.
	\begin{proof}
		The proof is similar to that of Theorem \ref{CtrDolbeault} but easier since the inclusion is already compatible with the dg-Lie structure.
	\end{proof}
\end{theorem}

\subsection{Matrix factorisations}
For a dg-algebra $A$ and $f \in A^0$ a central closed element, one may define the $\mathbb Z/2$-graded dg-category of matrix factorisations as the full subcategory of the dg-category of all modules over $(A, 0, f)$ (viewed as a curved $\mathbb Z/2$-graded dg-algebra) such that the two underlying $A$-modules are finitely generated projective over $A$. We note that all our definitions and results also work if we replace our grading group by $\mathbb Z/2$.

It is known that the Hochschild cohomology of the category of matrix factorisations of a smooth variety $X$ over $\mathbb C$ endowed with a Landau-Ginzburg potential $W$ is given by 
\[ \HH^*(\mathrm{MF}(X,W)) = \mathbb R \Gamma(\Lambda^* T_X, [W,-])\]
see \cite[Theorem 3.1]{MF} or \cite[Theorem 8.2.6]{Preygel}. We have analogous results in the smooth/analytic setting:
\begin{theorem}
	\begin{enumerate}
		\item	Let $f \in C^{\infty}(M)$ be a smooth function on a smooth manifold $M$. Then the continuous Hochschild cohomology of the $\mathbb Z/2$ clct cdg-algebra $(C^\infty(M),0,f)$ is given by
	\[ \HHcont(C^\infty(M),0,f) = H^* (\Lambda^* \Gamma(TM)[-1], d = \iota_{df}).\]
	\item Now let $A = \mathbb K[[x_1,...,x_n]]$ be the pseudo-compact algebra of formal power series in $n$ variables and $f \in A$. Then we have 
	\[ \HHcont(A,0,f) = H^* (\Lambda_{A}^* \text{Der}(A)[-1], d = \iota_{df})\]
	(note that $\text{Der}(A) \cong \mathbb K^n \otimes A$).
	\item Similarly, for the discrete polynomial algebra $P = \mathbb K[x_1,...,x_n]$ we obtain 
	\[ \HHcont(A,0,f) = H^* (\Lambda_{P}^* \text{Der}(P)[-1], d = \iota_{df})\]
	with $\text{Der}(P) \cong \mathbb K^n \otimes P$.
	
	\item Finally, for $X$ a complex manifold and $f \in \mathcal O_X(X)$ a global holomorphic function on $X$ we have
	\[ \HHcont(\mathcal A(X),\bar \partial, f) = H^*(\Lambda T^{1,0} X \otimes \mathcal A(X)[-1], \bar \partial + \iota_{\partial f} \otimes 1).\]
	\end{enumerate}
	\begin{proof}
		The proof of the first three results are similar (when $M = \mathbb R^n$) and the globalisation of the first statement can be carried out in complete analogy to Theorem \ref{HochschildFunctions}. We prove the second statement. The proof idea is very simple: We take the standard Koszul resolution of $A$ as a $A^{op} \hat \otimes A = \mathbb K[[x_1,...,x_n,y_1,...,y_n]]$-module and then twist this resolution to incorporate $f$. We thus consider:
		\[ 0 \to \Lambda^{n} \mathbb K^n \otimes \mathbb K[[x_1,...,x_n,y_1,...,y_n]] \to ... \to \Lambda^{0} \mathbb K^n \otimes \mathbb K[[x_1,...,x_n,y_1,...,y_n]] \to \mathbb K[[z_1,...,z_n]] \to 0 \]
		where the differential is given by
		\[ d(\alpha \otimes g) := \sum_{i = 1}^n \iota_{e^i} \alpha \otimes (x_i - y_i) g.\]
		This sequence is horizontally admissibly exact and provides a graded-free bimodule resolution of $A$ (when $f = 0$). To make it a bimodule resolution of $(A,0,f)$ we add differentials \[k : \Lambda^* \mathbb K^n  \otimes \mathbb K[[x_1,...,x_n,y_1,...,y_n]] \to \Lambda^{*+1} \mathbb K^n  \otimes \mathbb K[[x_1,...,x_n,y_1,...,y_n]] \]
		(note that they now increase the exterior degree). Of course $k$ is of degree $1 = -1$ in the $\mathbb Z/2$-graded setting. $k$ is constructed in the following manner. The curvature element of $A^{op} \hat \otimes A$ is given by $F(x,y) = f(y) - f(x)$. This implies that $F$ can be written as \begin{nalign}
			\label{ExpansionFormula}
		F(x,y) = \sum_{i = 1}^n (y_i - x_i) F_i(x,y).\end{nalign} Then we set
		\begin{nalign}\label{Expansion2} k(\alpha \otimes g) := \sum_{i = 1}^n e^i \wedge \alpha \otimes g F_i.
		\end{nalign}
		Then one finds that $k^2 = 0, (d\circ k + k \circ d)(\alpha \otimes g) = \alpha \otimes g F$ so that \\$(\Lambda^* \mathbb K^n \otimes \mathbb K[[x_1,...,x_n,y_1,...,y_n]],d+k)$ is indeed a module over $(A^{op} \hat \otimes A, 0, F)$. Clearly, it is graded-free and it is isomorphic (in the contraderived category) to $\mathbb K[[z_1,...,z_n]]$ by Lemma \ref{lemmaGenResolution}. The statement follows by applying $\Hom(-,A)$ and noting that $\Delta^* F_i(x) = F_i(x,x) = \partial_{x_i} F(x)$.
		
		We now turn to the fourth statement. As before, the proof can be reduced to a local computation. In this case, it is useful to consider (open) polydiscs $\Delta \subset M$ and the corresponding opens $\Delta \times \Delta \subset M$. By the proof of Proposition \ref{LocalDolbeault}, we obtain a (Koszul-type) resolution $R^{std,-*}$ of $\mathcal A(\Delta)$ as a module over $\mathcal A(\Delta \times \Delta)$. The set $\Delta \times \Delta$ is convex, and so Taylor expansion yields functions $F_i$ on $\Delta \times \Delta$ as in formula \ref{ExpansionFormula} and we can twist $R^{std,-*}$ as in \ref{Expansion2} to make it into a $\mathbb Z_2$ resolution of $\mathcal A(\Delta)$ over $(\mathcal A(\Delta \times \Delta),\bar \partial, F)$ in the sense of Lemma \ref{lemmaGenResolution}. Analogously to the proof of \ref{LocalDolbeault}, we find that
		\[ \text{Hom}^*_{(\mathcal A(\Delta \times \Delta),\bar \partial, F)}(R^{std,-*}_{tw}, \mathcal A(\Delta)) \simeq  (\Lambda T^{1,0} \Delta \otimes \mathcal A(\Delta)[-1], \bar \partial + \iota_{\partial f} \otimes 1).\]
	\end{proof}
\end{theorem}

\subsection{Outlook: Formality for Hochschild complexes of enhanced derived categories}
Let $A$ be a clct (c)dg-algebra. The dg-category $\text{D}^{\text{ctr},b}(A)$ of cdg-modules that are finitely generated and projective as graded modules is naturally enriched in clct dg-vector spaces. In the follow-up paper, we define the continuous Hochschild complex for such categories and show an invariance result of the $B_\infty$ type of the Hochschild complex under so called contraderived equivalences. This result in particular implies that the continuous Hochschild-complex of $\text{D}^{\text{ctr},b}(A)$ is $B_\infty$-equivalent to that of $A$. We thus in particular obtain:
\begin{theorem}
	\label{FormalityEnhancement}
	Let $X$ be a compact complex manifold. Then the bounded derived category $\D{b}{coh}(X)$ has a natural clct dg-enhancement $\D{ctr,b}{}(\mathcal A(X))$ where $\mathcal A(X)$ is the Dolbeault algebra. Moreover, there is an $L_\infty$-quasi-isomorphism
	\[ \HC(\D{ctr,b}{}(\mathcal A(X))) \simeq \Lambda T^{1,0} X \otimes \mathcal A(X)[-1].\]
	The right-hand side is isomorphic to the deformation complex of $X$ as a generalized complex manifold.
	If $X$ is not compact, then the above statement holds if we replace $\D{b}{coh}(X)$ by the category of globally bounded perfect complexes on $X$.
\end{theorem}
Of course, similar statements hold for the other examples discussed in the previous sections. In particular, the continuous Hochschild complex of the de-Rham clct enhancement of the category of $\infty$-local systems on a smooth manifold $M$ is dg-Lie quasi-isomorphic to the (shifted) de-Rham complex (with trivial bracket).

\appendix
\section{Background on complete locally convex vector spaces}
We will collect some facts about complete locally convex spaces over $\mathbb K = \mathbb R$ or $\mathbb K = \mathbb C$ here.

\begin{Sdef1}
	\label{projTensorProduct}
		Let $M$ and $N$ be complete locally convex topological vector spaces (clctvs) endowed with families of seminorms $|| \cdot ||_i$ and $|| \cdot ||_j$. Their \emph{uncompleted projective tensor product $M \otimes_{\pi} N$} is the algebraic tensor product of $M$ and $N$ together with the family of seminorms 
		\[ ||f||_{i,j} = \inf\{\sum ||f_a||_i ||f_b||_j \text{ } | \text{  } f = \sum f_a \otimes f_b\}\]
		where $(i,j) \in I \times J$. Its completion will be called the \emph{projective tensor product} or the \emph{completed tensor product} and we denote it by
		\[ M \hat \otimes N := \widehat{M \otimes_{\pi} N}\]
		which is again a clctvs.
\end{Sdef1}
\begin{Slemma}
	\label{ProjTensor}
	Let $M_{k \in K}$ be a family of clctvs and $N$ another clctvs.
	Then there is a canonical isomorphism
	\[ (\prod_{k \in K} M_k) \hat \otimes N = \prod_{k \in K} (M_k \hat \otimes N).\]
	\begin{proof}
		See \cite[Section 15.4]{Jarchow}
	\end{proof}
\end{Slemma}
\begin{Slemma}
	\label{InductiveLimits}
	Let $M_{m \in \mathbb N}$ and $N_{n \in \mathbb N}$ be countable inductive systems of Banach spaces.
	Then there is a canonical isomorphism
	\[ (\lim\limits_{\longrightarrow} M_m) \hat \otimes (\lim\limits_{\longrightarrow} N_n) =\lim\limits_{\longrightarrow}  (M_m \hat \otimes N_n).\]
	\begin{proof}
		See \cite[Corollary 15.5.4]{Jarchow}.
	\end{proof}
\end{Slemma}
\begin{Slemma}
	\label{LemmaProdPlus}
	Let $M_{k \in K}$ be a sequence of clctvs and $B$ a Banach space (or more generally, if $B$ admits a Hausdorff-seminorm). Then
	\begin{nalign}
		\label{BanachCocompact}
		\Hom_{\mathbb C}(\prod_{k \in K} M_k, B) = \oplus_{k \in K} \Hom_{\mathbb C}(M_k, B) 
	\end{nalign}
	\begin{proof}
		We use the standard description of seminorms on $\prod_{k \in K} M_k$. Given a map $f : \prod_{k \in K} M_k \to B$, the one and only norm on $B$ is bounded by finitely many seminorms on $\prod_{k \in K} M_k$. Each of these seminorms is associated to some factor $M_k$ so we get a finite collection $K' \subseteq K$ so that $f$ factors through $\prod_{k \in K} M_k \to \prod_{k \in K'} M_k$. 
	\end{proof}
\end{Slemma}
\begin{Sdef1}
	\label{definitionLocalisation}
	Let $N$ be a clctvs. Given a seminorm $||\cdot||_j$ on $N$ we consider the Banach space $N_j$ which is the completion of 
	\[ N_j^{pre} = (N / \{v \in N, ||v||_{j} = 0\}, ||\cdot||_j)\]
	and we denote by $p_j$ the composition $p_j : N \to N_j^{pre} \to N_j$.
	
	A reduced inverse limit of complete locally convex spaces is an inverse limit ${\lim\limits_{\longleftarrow} N^j}$ for which all maps ${\lim\limits_{\longleftarrow} N^j} \to N^i$ have dense range ($\Leftrightarrow$ the cokernel is trivial). 
\end{Sdef1}
\begin{Slemma}
	Every clctvs is a reduced inverse limit of Banach spaces. 
	\begin{proof}
		See \cite[Theorem 6.8.5]{Jarchow}. 
	\end{proof}
\end{Slemma}
\begin{Scorollary}
	\label{CorProdHom}
	Let $M_{k \in K}$ be a family of clctvs and $N$ another clctvs. Write $N$ as a reduced inverse limit of Banach spaces $N = {\lim\limits_{\longleftarrow} N_j}$. Then $\Hom_{\mathbb C}(\prod_{k \in K} M_k, N)$ is the subspace
	\[ \Hom_{\mathbb C}(\prod_{k \in K} M_k, N) \subseteq \prod_{k \in K} \Hom_{\mathbb C}(M_k, N)\]
	consisting of those families $(f_k : M_k \to N)_{k \in K}$ with the property that the families $(p_j \circ M_k \to N_j)_{k \in K}$ lie in $\oplus_{k \in K} \Hom_{\mathbb C}(M_k, N_j)$ for every $j \in J$. In formulas:
	\begin{nalign}  \Hom_{\mathbb C}(\prod_{k \in K} M_k, N) = \lim_{\stackrel{\longleftarrow}{j \in J}} \bigoplus_{k \in K} \Hom_{\mathbb C}(M_k, N_j).\end{nalign}
\end{Scorollary}

\begin{Scorollary}
	\label{ProdCokernels}
	Arbitrary products commute with cokernels in the category of complete locally convex vector spaces.

	\begin{proof}
		Let $(M_k \to N_k)_{k \in K}$ be a collection of maps between complete locally convex vector spaces. As every complete locally convex space is an inverse limit of Banach spaces, it suffices to show that for each Banach space $B$ there is a canonical isomorphism
		\begin{nalign} \label{EqCokerKer} \Hom(\prod_{k \in K} \text{coker}(M_k \to N_k), B) =  \Hom(\text{coker}(\prod_{k \in K} M_k \to \prod_{k \in K} N_k), B).
		\end{nalign}
		By Lemma \ref{LemmaProdPlus},
		the left hand side is \[\oplus_{k \in K} \Hom( \text{coker}(M_k \to N_k), B) = \oplus_{k \in K} \text{ker} (\Hom( N_k, B) \to \Hom( M_k, B)).\] Now as the abelian category of discrete vector spaces fulfills AB4, the latter expression is equal to
		\[ \text{ker} ( \oplus_{k \in K} \Hom( N_k, B) \to  \oplus_{k \in K} \Hom( M_k, B))\]
		which equals the right-hand side of equation \ref{EqCokerKer}.
	\end{proof}
\end{Scorollary}

\begin{Slemma}
	\label{DiscreteLocallyConvex}
	Let $V$ be a (discrete) vector space. There exists the finest topology on $V$ making $V$ into a complete locally convex vector space. Every seminorm on $V$ is then continuous. It can also be described by choosing a (Hamel) basis $B$ for $V$. Then $V = \oplus_B \mathbb K$ where the sum is taken in the category of complete locally convex vector spaces. In this way we obtain a fully-faithful functor which is left-adjoint to the forgetful functor from complete locally convex spaces to (discrete) vector spaces.
	\begin{proof}
		See \cite[p. 111]{Jarchow}.
	\end{proof}
\end{Slemma}

\bibliographystyle{alphaurl}
\bibliography{pmetemplate03}

\end{document}